\documentclass[12pt]{amsart}
\usepackage{amscd,amsfonts,amsmath,amsxtra,amssymb}

\usepackage[all]{xy}

\newtheorem{sub}{}[section]
\newtheorem{subsub}{}[sub]

\newcommand{\R}{{\mathbb R}}
\newcommand{\C}{{\mathbb C}}
\newcommand{\Q}{{\mathbb Q}}
\renewcommand{\P}{{\mathbb P}}

\def\T{{\mathbb T}}

\newcommand{\ki}{{\mathcal I}}

\newcommand{\km}{{\mathcal M}}

\newcommand{\ko}{{\mathcal O}}

\newcommand{\kq}{{\mathcal Q}}

\def\lra{\longrightarrow}
\def\sigg{\mathop{\hbox{$\displaystyle\sum$}}\limits}

\def\ot{\otimes}

\def\and{\quad\text{and}\quad}

\def\dsp{\displaystyle}

\font\tte=cmbsy10

\def\ov#1{\overline{#1}}

\def\Ext{\mathop{\rm Ext}\nolimits}
\def\Aut{\mathop{\rm Aut}\nolimits}
\def\End{\mathop{\rm End}\nolimits}
\def\Cong{\mathop{\rm Cong}\nolimits}
\def\Isom{\mathop{\rm Isom}\nolimits}

\def\lra{\longrightarrow}

\def\sigg{\mathop{\hbox{$\displaystyle\sum$}}\limits}

\def\timex{\setbox250=\hbox{$\times$}\hskip 5pt\Box\hskip -0.75em
{\raise 1.5pt\vbox{\box250}}\hskip 5pt}
\def\paragra{{\tte \char120}}
\def\para{\paragra~\hskip -2pt}

\def\hfl#1#2{\smash{\mathop{\hbox to 12mm{\rightarrowfill}}
\limits^{\scriptstyle#1}_{\scriptstyle#2}}}

\def\pline#1{<\hskip-3.5pt#1\hskip-3.5pt>}

\def\m#1{{\hbox{$#1$}}}
\def\ot{\otimes}
\def\og{\leavevmode\raise.3ex\hbox{$\scriptscriptstyle\langle\!\langle$}}
\def\fg{\leavevmode\raise.3ex\hbox{$\scriptscriptstyle\,\rangle\!\rangle$}}

\begin{document}

\title[{\tiny G\'eom\'etrie orthogonale non sym\'etrique}]
{G\'eom\'etrie orthogonale non sym\'etrique et congruences quadratiques}
\def\refname{R\'ef\'erences}
\def\contentsname{Sommaire}
\def\proofname{D\'emonstration}
\def\abstractname{R\'esum\'e}

\author[{\tiny J.M. Dr\'{e}zet}]{Jean--Marc Dr\'{e}zet}
\address{
Institut de Math\'{e}matiques\\
UMR 7586 du CNRS, Aile 45-55, 5$^e$ \'{e}tage\\ 
2, place Jussieu\newline   
F-75251 Paris Cedex 05, France}
\email{drezet@math.jussieu.fr}

\maketitle
\tableofcontents
\section{Introduction}

Le but de cet article est l'\'etude de certaines structures g\'eom\'etriques
rencontr\'ees dans la th\'eorie des fibr\'es exceptionnels sur les espaces
projectifs (cf. \cite{dr3_0}, \cite{dr3}). Le point de d\'epart est la donn\'ee
d'un $\C$-espace vectoriel $W$
de dimension finie \m{n+1}, \m{n\geq 3}, muni d'une forme bilin\'eaire
non d\'eg\'en\'er\'ee. Dans l'\'etude des fibr\'es exceptionnels sur \m{\P_n},
on consid\`ere \ \m{W={\rm K}(\P_n)\ot\C} \ (\m{{\rm K}(\P_n)} d\'esignant
l'{\em{anneau de Grothendieck}} des classes de faisceaux coh\'erents sur
\m{\P_n}), muni de la forme bilin\'eaire d\'efinie par
\[\chi(\lbrack E\rbrack, \lbrack F\rbrack) \ = \ \sigg_{i=0}^n\dim(\Ext^i(E,F))
\]
pour tous faisceaux coh\'erents $E$, $F$ sur \m{\P_n}. L'\'etude de cette forme
bilin\'eaire et des g\'en\'e-\break ralisations ont \'et\'e effectu\'ees dans
\cite{go3}, \cite{ru3}. Le cas de \m{\P_3} est trait\'e avec beaucoup de
d\'etails dans \cite{no1}, \cite{no2}. Les articles pr\'ec\'edents \'etaient
plut\^ot consacr\'es aux propri\'et\'es arithm\'etiques de $\chi$ sur 
\m{{\rm K}(\P_n)}. On s'int\'eresse ici uniquement \`a des formes bilin\'eaires
complexes.

\medskip

Soit \ \m{\phi : W\lra W^*} \ un isomorphisme non antisym\'etrique. 
On en d\'eduit la forme
bilin\'eaire non d\'eg\'en\'er\'ee d\'efinie par
\[(u,v) \ = \ u.\phi(v) .\]
On note \ \m{Q(\phi)\subset\P(W)} \ la quadrique des points isotropes. 
Si $X$ est un sous-ensemble de \m{\P(W)}, on note \m{X^\bot} (resp. 
\m{{^\bot}X}) le sous-espace
lin\'eaire de \m{\P(W)} constitu\'e des points $y$ tels que \ \m{(x,y)=0} 
(resp. \m{(y,x)=0}) pour tout $x$ dans $X$.
Soit \ \m{\T(\phi)={^t}\phi^{-1}\circ\phi\in GL(W)}, appel\'e ici la {\em
translation} associ\'ee \`a $\phi$ \ (dans \cite{go3}, \m{\T(\phi)^{-1}} est
appel\'e l'{\em{op\'erateur fondamental}} associ\'e \`a la forme bilin\'eaire
pr\'ec\'edente). On rappelle au chapitre 2 quelques r\'esultats classiques, et
notamment celui-ci : la \m{GL(W)}-orbite de $\phi$ est enti\`erement
d\'etermin\'ee par celle de \m{\T(\phi)} \ (cf. \cite{go3}, et les ouvrages
\cite{Ma}, \cite{H_P} cit\'es dans cet article). On donne au \para 2.2.2 la
liste des orbites des \m{\T(\phi)} et celle des $\phi$ correspondants dans le
cas o\`u \ \m{n=3}.

\medskip

Soit $V$ un $\C$-espace vectoriel de dimension $n$. On suppose que $W$ est un
sous-espace vectoriel de \m{H^0(\P(V),\ko(2))}. On appelle {\em congruence
quadratique} un morphisme rationnel
\[\sigma : \P(V)\lra\P(W)\]
tel que \ \m{\sigma^*(\ko(1))=\ko(2)}, que \ \m{\sigma^*:W^*\lra 
H^0(\P(V),\ko(2))}
\ induise un isomorphisme \ \m{\sigma^*:W^*\simeq W} \ et que pour tout point 
$P$ de \m{\P(V)} o\`u $\sigma$ est d\'efini on ait \m{P\in\sigma(P)} \ 
(\m{\sigma(P)} \'etant vu comme une quadrique de \m{\P(V)}). Soit \
\m{\phi=(\sigma^*)^{-1}}. On munit alors $W$ de la forme bilin\'eaire non
d\'eg\'en\'er\'ee induite par $\phi$. L'\'etude des congruences quadratiques
dans le cas g\'en\'eral est faite dans le chapitre 3. Dans le chapitre 4 on
s'int\'eresse au cas \m{n=3} (on parle alors de congruences quadratiques {\em
planes}). 

Le cas le plus int\'eressant est
celui o\`u \m{\P(W)} est l'espace des quadriques de \m{\P(V)} contenant une
quadrique fixe \m{C(\sigma)} d'un hyperplan \m{H(\sigma)} de \m{\P(V)}. 
On dit dans ce cas que $\sigma$
est une congruence quadratique {\em normale}, et on montre qu'on a alors \
\m{rg(C)=rg(Q(\phi))-2}, que $\sigma$ induit un isomorphisme birationnel \
\ \m{\P(V)\simeq\ov{\sigma(\P(V))}} \ et que \ \m{\ov{\sigma(\P(V))}=Q(\phi)}.
On note
\[\T(\sigma) \ = \ \sigma^{-1}\circ\T(\phi)\circ\sigma,\]
qui est un automorphisme birationnel de \m{\P(V)}, qu'on appelle la {\em
translation} associ\'ee \`a $\sigma$. 

La restriction de $\sigma$ \`a \m{H(\sigma)} est constante. C'est une quadrique
d\'eg\'en\'er\'ee dont l'une des composantes est \m{H(\sigma)}. On note
\m{L(\sigma)} l'autre composante. On montre que \m{\T(\sigma)} est lin\'eaire si
et seulement si on a \ \m{L(\sigma)=H(\sigma)}.

\medskip

On donne au \para 3.5 une construction g\'eom\'etrique des congruences
quadratiques normales. On part d'un espace vectoriel abstrait $W$ de dimension
\m{n+1} et d'un
isomorphisme non antisym\'etrique \ \m{\phi : W\simeq W^*} tel que \m{Q(\phi)}
soit non d\'eg\'en\'er\'ee. Soient $O$ un point
lisse de \m{Q(\phi)}, \m{T_O} l'hyperplan tangent \`a \m{Q(\phi)} en $O$ et
\[\pi : \P(W)\backslash\lbrace O\rbrace\lra\P_{n-1}\]
la projection de centre $O$. Soit
\[C \ = \ \pi(T_O\cap Q(\phi))\]
qui est une quadrique de l'hyperplan \m{\pi(T_O)} de \m{\P_{n-1}}. Soit \m{W_0}
le sous-espace vectoriel de \m{H^0(\P_{n-1},\ko(2))} constitu\'e des \'equations
de quadriques contenant $C$. On d\'efinit alors un isomorphisme \
\m{\P(W)\simeq\P(W_0)} \ en associant \`a $P$ l'image par $\pi$ de
\m{Q(\phi)\cap P^\bot}. Si $Q$ est un point g\'en\'eral de \m{\P_{n-1}}, on note
\m{\eta(Q)} l'unique point de \m{Q(\phi)\cap\pi^{-1}(Q)}. Alors
\[
\begin{array}{cccc}
\sigma : & \P_{n-1} & \lra        & \P(W_0)        \\
         & Q        & \longmapsto & \pi(Q(\phi)\cap{^\bot\eta(Q)})
\end{array}
\]
est une congruence quadratique. On montre (th\'eor\`eme \ref{theo3}) que toute
congruence quadratique normale
peut s'obtenir de cette fa\c con (modulo l'action de \m{PGL(n)}).

\medskip

Dans le chapitre 4 on suppose que \ \m{n=3}. Dans ce cas $C$ est constitu\'e de
deux points distincts de \m{\P_2}, ou c'est un point double d'une droite de
\m{\P_2}. Soient \m{G_C} le sous-groupe de \m{GL(3)} cons-\break titu\'e des 
isomorphismes laissant $C$ invariant et \m{\Cong(C)} l'ensemble des congruences
quadratiques \ \m{\sigma:\P_2\lra\P(W)}, o\`u $W$ est l'espace des \'equations
de coniques de \m{\P_2} contenant $C$. Alors l'application
\[ \begin{array}{ccc} \Cong(C) & \lra & \Isom(W,W^*) \\
\sigma  & \longmapsto & (\sigma^*)^{-1} \\ \end{array}\]
est compatible avec un morphisme de groupes \ \m{G_C\lra GL(W)}. On montre 
(th\'eor\`eme \ref{theo4}) que l'application quotient 
\[G_C/G_C\lra\Isom(W,W^*)/GL(W)\]
a des fibres finies. On donne au chapitre 6 la liste des \m{G_C}-orbites de
\m{\Cong(C)} correspondant aux \m{GL(W)}-orbites des \m{\T(\phi)}
correspondantes.

\medskip

Dans le chapitre 5 on montre (th\'eor\`eme \ref{th0b}) que pour toute congruence
quadratique normale \ \m{\sigma:\P(V)=\P_{n-1}\lra\P(W)} \ telle que l'image de
$\sigma$ contienne au moins une quadrique de rang maximal de \m{\P(W)} 
il existe un morphisme
rationnel \ \m{R:\P(V)\lra PGL(V)} \ et une quadrique \m{C_0} de \m{\P(V)}
tels que pour un point g\'en\'eral $P$ de \m{\P(V)} on ait \ 
\m{\sigma(P)=R^{-1}(C_0)}. Cela d\'ecoule du fait qu'il existe des sections
rationnelles de la restriction \`a \m{\P(W)} de la quadrique universelle \
\m{\kq\subset\P(S^2V^*)\times\P(V)} (prop. \ref{prop8}). Dans le cas o\`u \
\m{n=3}, on montre que si \ \m{L(\sigma)=H(\sigma)}, on peut d\'efinir $R$ par
des formes de degr\'e \m{\leq 2}.

\vskip 2.5cm

\section{Isomorphismes, quadriques et translations associ\'ees}

\begin{sub}\label{defgen}\bf D\'efinitions

\rm Soit $W$ un \m{\C}-espace vectoriel de dimension finie \ \m{n+1\geq 4}. Soit
\[
\begin{array}{cccc}
\T : & \Isom(W,W^*) & \lra & GL(W)\\
    & \quad\quad\quad\quad\phi         & \longmapsto & {}^t\phi^{-1}\circ\phi
\end{array}
\]
Dans \cite{dr3} on appelle \m{\T(\phi)} la {\em translation} associ\'ee \`a 
$\phi$. Dans \cite{go3} on appelle \m{\T(\phi)^{-1}} 
l'{\em op\'erateur canonique}.
Le groupe \m{GL(W)} agit sur lui-m\^eme \`a droite par conjugaison et sur 
\m{\Isom(W,W^*)} par :
\[
\phi.g \ = \ ^tg\circ\phi\circ g
\]
pour \ \m{g\in GL(W), \phi\in\Isom(W,W^*)}, et $\T$ est un $GL(W)$-morphisme.

Si \ \m{\phi\in\Isom(W,W^*)}, on note 
\[
\phi_S=\frac{1}{2}(\phi+{^t\phi}), \ \ \ \phi_A=\frac{1}{2}(\phi-{^t\phi})
\]
qu'on appelle respectivement la {\em partie sym\'etrique} et la {\em partie
antisym\'etrique} de $\phi$. On utilise les m\^emes notations pour les
\'el\'ements de \m{\Isom(W^*,W)}.

On notera souvent de la m\^eme fa\c con un point de \m{\P_n} et un point de $W$
au dessus du pr\'ec\'edent, ce qui donne par exemple un sens \`a la 
d\'efinition suivante : si $\phi$ est non antisym\'etrique, on note 
\m{Q(\phi)} la
quadrique de \ \m{\P_n=\P(W)} \ d'\'equation \ \m{x.\phi(x)=0}. 

On note \m{X^\bot} (resp. \m{^\bot X}) le sous-espace lin\'eaire de \m{\P_n} 
constitu\'e des points \m{y} tels que \ \m{x\phi(y)=0} (resp. \ 
\m{y\phi(x)=0}) pour tout point $x$ de $X$.

\vskip 1cm

\begin{subsub}\label{lemm8}{\bf Lemme : } Soit $P$ un point lisse de 
\m{Q(\phi)}. Alors les trois propri\'et\'es suivantes sont \'equivalentes :
\begin{enumerate}
\item On a \ $\phi(P)={^t\phi}(P)$ \ dans $\P(W^*)$.
\item On a \ $P^\bot = {^\bot}P$.
\item On a \ $\T(\phi)(P) = P$ \ dans $\P(W)$.
\item L'espace tangent de \m{Q(\phi)} en $P$ est $P^\bot$. 
\end{enumerate}
\end{subsub}
\end{sub}

Imm\'ediat

\vskip 1.5cm

\begin{sub}\label{rap1}{\bf \'Etude de l'action de \m{GL(W)}}\rm

La \m{GL(W)}-orbite d'un 
\'el\'ement $\phi$ de \m{\Isom(W,W^*)} est enti\`erement d\'etermin\'ee par la
\m{GL(W)}-orbite de \m{\T(\phi)} :

\vskip 1cm

\begin{subsub}\label{theo1}{\bf Proposition : } 
Pour tout \ \m{\phi\in\Isom(W,W^*)}, on a
\[
\T^{-1}(\T(\phi)) \ \subset \ GL(W).\phi.
\]
\end{subsub}

(Proposition 3.1.1 de \cite{go3}).

\bigskip


\vskip 1cm

\begin{subsub}\label{liste_0}Le cas o\`u \ \m{\dim(W)=4}\end{subsub}

On donne ci-dessous la liste des orbites des \m{\T(\phi)} et la forme des $\phi$
correspondants dans le cas o\`u \ \m{\dim(W)=4}. 

\medskip

\noindent Cas 1.1 : 
\[
\T(\phi)=\left(\begin{array}{cccc}
\lambda & 0 & 0 & 0 \\ 0 & \frac{1}{\lambda} & 0 & 0\\
0 & 0 & \mu & 0 \\ 0 & 0 & 0 & \frac{1}{\mu}
\end{array}\right),\quad\quad\phi=\left(\begin{array}{cccc}
0 & X & 0 & 0 \\ \lambda X & 0 & 0 & 0 \\ 0 & 0 & 0 & Y \\ 0 & 0 & \mu Y & 0
\end{array}\right)
\]
avec  \  \m{\lambda^2\not=1,\ \mu\not=1,\ \lambda\mu\not=1,\ \lambda\not=\mu},
\ \m{XY\not=0}. On a \ \m{rg(Q(\phi))=4}.

\medskip

\noindent Cas 1.2 :
\[
\T(\phi)=\left(\begin{array}{cccc}
\lambda & 0 & 0 & 0 \\ 0 & \frac{1}{\lambda} & 0 & 0 \\
0 & 0 & \lambda & 0 \\ 0 & 0 & 0 & \frac{1}{\lambda}
\end{array}\right),\quad\quad\phi=\left(\begin{array}{cccc}
0 & X & 0 & T \\ \lambda X & 0 & \lambda Z & 0 \\ 
0 & Z & 0 & Y \\ \lambda T & 0 & \mu \lambda Y & 0
\end{array}\right)
\]
avec  \  \m{\lambda^2\not=1},\ \m{XY-ZT\not=0}. On a \ \m{rg(Q(\phi))=4}.

\medskip

\noindent Cas 1.3 :
\[
\T(\phi)=I_W, \quad \phi=\ \text{matrice sym\'etrique},
\]
On a \ \m{rg(Q(\phi))=4}.

\medskip

\noindent Cas 1.4 :
\[
\T(\phi)=-I_W, \quad \phi=\ \text{matrice antisym\'etrique},
\]

\medskip

\noindent Cas 1.5 :
\[
\T(\phi)=\left(\begin{array}{cccc}
1 & 0 & 0 & 0 \\ 0 & 1 & 0 & 0\\
0 & 0 & \lambda & 0 \\ 0 & 0 & 0 & \frac{1}{\lambda}
\end{array}\right),\quad\quad\phi=\left(\begin{array}{cccc}
\alpha & X & 0 & 0 \\ X & \beta & 0 & 0 \\ 
0 & 0 & 0 & Y \\ 0 & 0 & \lambda Y & 0
\end{array}\right)
\]
avec  \  \m{\lambda\not=1},\ \m{(\alpha\beta-X^2)Y\not=0}. 
On a \ \m{rg(Q(\phi))=4} \ si \ \m{\lambda\not=-1} \ et sinon \ 
\m{rg(Q(\phi))=2}.

\medskip

\noindent Cas 2.1 :
\[
\T(\phi)=\left(\begin{array}{cccc}
-1 & 1 & 0 & 0 \\ 0 & -1 & 1 & 0 \\ 0 & 0 & -1 & 1 \\ 0 & 0 & 0 & -1
\end{array}\right),\quad\quad\phi=\left(\begin{array}{cccc}
0 & 0 & 0 & \alpha \\ 0 & 0 & -\alpha & -\frac{\alpha}{2} \\ 
0 & \alpha & -\frac{\alpha}{2} & \beta \\ 
-\alpha & \frac{3}{2}\alpha & -\frac{\alpha}{2}-\beta & \frac{\beta}{2}
\end{array}\right)
\]
avec \ \m{\alpha\not=0}. On a \ \m{rg(Q(\phi))=3}.

\medskip

\noindent Cas 2.2 :
\[
\T(\phi)=\left(\begin{array}{cccc}
-1 & 1 & 0 & 0 \\ 0 & -1 & 0 & 0 \\ 0 & 0 & -1 & 1 \\ 0 & 0 & 0 & -1
\end{array}\right),\quad\quad\phi=\left(\begin{array}{cccc}
0 & 2\alpha & 0 & u+v \\ -2\alpha & \alpha & -u-v & v \\ 
0 & u+v & 0 & 2\beta \\ -u-v & u & -2\beta & \beta
\end{array}\right)
\]
avec \ \m{4\alpha\beta-(u+v)^2\not=0}. On a \ \m{rg(Q(\phi))=2}.

\medskip

\noindent Cas 2.3 :
\[
\T(\phi)=\left(\begin{array}{cccc}
\lambda & 1 & 0 & 0 \\ 0 & \lambda & 0 & 0 \\
0 & 0 & \frac{1}{\lambda} & 1 \\ 0 & 0 & 0 & \frac{1}{\lambda}
\end{array}\right),\quad\quad\phi=\left(\begin{array}{cccc}
0 & 0 & 0 & -\lambda^2\alpha \\ 0 & 0 & \alpha & \beta \\ 
0 & \lambda\alpha & 0 & 0 \\ 
-\lambda^3\alpha & -\lambda^2\alpha+\lambda\beta & 0 & 0
\end{array}\right)
\]
avec \ \m{\lambda^2\not=1}, \m{\alpha\not=0}. On a \ \m{rg(Q(\phi))=4}.

\medskip

\noindent Cas 2.4 :
\[
\T(\phi)=\left(\begin{array}{cccc}
1 & 1 & 0 & 0 \\ 0 & 1 & 0 & 0 \\ 0 & 0 & 1 & 1 \\ 0 & 0 & 0 & 1
\end{array}\right),\quad\quad\phi=\left(\begin{array}{cccc}
0 & 0 & 0 & \alpha \\ 0 & \epsilon & -\alpha & \beta-\alpha \\ 
0 & -\alpha & 0 & 0 \\ \alpha & \beta & 0 & \gamma
\end{array}\right)
\]
avec \ \m{\alpha\not=0}. On a \ \m{rg(Q(\phi))=4}.

\medskip

\noindent Cas 2.5 :
\[
\T(\phi)=\left(\begin{array}{cccc}
1 & 0 & 0 & 0 \\ 0 & 1 & 1 & 0 \\ 0 & 0 & 1 & 1 \\ 0 & 0 & 0 & 1
\end{array}\right),\quad\quad\phi=\left(\begin{array}{cccc}
\alpha & 0 & 0 & \beta \\ 0 & 0 & 0 & \gamma \\ 
0 & 0 & -\gamma & 0 \\ \beta & \gamma & \gamma & \epsilon
\end{array}\right)
\]
avec \ \m{\alpha\gamma\not=0}. On a \ \m{rg(Q(\phi))=4}.

\medskip

\noindent Cas 2.6 :
\[
\T(\phi)=\left(\begin{array}{cccc}
\lambda & 0 & 0 & 0 \\ 0 & \frac{1}{\lambda} & 0 & 0 \\
0 & 0 & -1 & 1 \\ 0 & 0 & 0 & -1
\end{array}\right),\quad\quad\phi=\left(\begin{array}{cccc}
0 & -\alpha & 0 & 0 \\ \lambda\alpha & 0 & 0 & 0 \\ 
0 & 0 & 0 & -\beta \\ 0 & 0 & \beta & -\frac{\beta}{2}
\end{array}\right)
\]
avec \ \m{\lambda^2\not=1}, \m{\alpha\beta\not=0}. On a \ \m{rg(Q(\phi))=3}.

\medskip

\noindent Cas 2.7 :
\[
\T(\phi)=\left(\begin{array}{cccc}
-1 & 0 & 0 & 0 \\ 0 & -1 & 0 & 0 \\ 0 & 0 & -1 & 1 \\ 0 & 0 & 0 & -1
\end{array}\right),\quad\quad\phi=\left(\begin{array}{cccc}
0 & -\alpha & 0 & 0 \\ \alpha & 0 & 0 & \gamma \\ 
0 & 0 & 0 & \beta \\ 0 & -\gamma & -\beta & \frac{\beta}{2}
\end{array}\right)
\]
avec \ \m{\alpha\beta\not=0}. On a \ \m{rg(Q(\phi))=1}.

\medskip

\noindent Cas 2.8 :
\[
\T(\phi)=\left(\begin{array}{cccc}
1 & 0 & 0 & 0 \\ 0 & 1 & 0 & 0 \\ 0 & 0 & -1 & 1 \\ 0 & 0 & 0 & -1
\end{array}\right),\quad\quad\phi=\left(\begin{array}{cccc}
a & c & 0 & 0 \\ c & d & 0 & 0 \\ 
0 & 0 & 0 & \alpha \\ 0 & 0 & -\alpha & \frac{\alpha}{2}
\end{array}\right)
\]
avec \ \m{\alpha(ad-c^2)\not=0}. On a \ \m{rg(Q(\phi))=3}.

\medskip

\noindent Les cas 1.x sont ceux o\`u \m{\T(\phi)} est diagonalisable 
les cas 2.x sont ceux o\`u elle ne l'est pas. Dans le cas g\'en\'eral les
r\'esultats de \cite{go3} permettent de d\'eterminer toutes les
\m{GL(W)}-orbites des $\phi$ et \m{\T(\phi)} correspondants.
\end{sub}

\newpage

\section{Congruences quadratiques (cas g\'en\'eral)}\label{congdef}

\begin{sub}\label{cong} \bf D\'efinition\rm

Soient \m{n\geq 3} un entier,
$V$ un $\C$-espace vectoriel de dimension \m{n}, et \ \m{\P_{n-1}=\P(V)}.
Soit $W_0$ un sous-espace vectoriel de dimension n+1 de \ 
\m{H^0(\P_{n-1},\ko(2))}. 
On appelle {\em congruence quadratique} un morphisme rationnel
\[
\sigma : \ \P_{n-1}\lra \P(W_0)=\P_n
\] 
tel que :
\begin{enumerate}
\item $\sigma^*(\ko(1)) = \ko(2)$ .
\item $\sigma^* : W_0^*\lra H^0(\P_{n-1}, \ko(2))$ \ induit un isomorphisme \ 
$W_0^*\simeq W_0$ \ aussi not\'e $\sigma^*$.
\item Pour tout point $P$ de $\P_2$ o\`u $\sigma$ est d\'efini, on a
\ $P\in\sigma(P)$, $\sigma(P)$ \'etant vu comme une quadrique de $\P_{n-1}$.  
\end{enumerate}
\medskip

Si de plus $\sigma$ induit un isomorphisme birationnel \ \m{\P_{n-1}\simeq
\ov{\sigma(\P_{n-1})}}, et si \m{\phi_S} est de rang au moins $n$
on dit que $\sigma$ est {\em r\'eguli\`ere}.
Si \ \m{n=3} \ on dit que $\sigma$ est une {\em congruence quadratique plane}.

Soit \ \m{\phi\in\Isom(W,W^*)}. On dit que $\sigma$ est une {\em r\'ealisation
g\'eom\'etrique} de $\phi$ {\em par des quadriques de \m{W_0}} s'il existe un 
isomorphisme \ \m{\eta : W\simeq W_0} \ tel que $\phi$ soit la compos\'ee
\[
\begin{array}{lclclcl}
W & \xrightarrow{\ \eta \ } & W_0 & \xrightarrow{(\sigma^*)^{-1}} & W_0^* &
\xrightarrow{\ ^t\eta \ } & W^*
\end{array}
\]
\end{sub}
On identifie dans ce cas $W$ et $W_0$, $\sigma^*$ et $\phi^{-1}$.

\vskip 1cm

\begin{subsub}{\bf Lemme : }\label{lemm1} 1 - Le morphisme $\sigma$ est \`a valeurs 
dans $Q(\phi)$, et pour tous points $P$,$Q$ de $\P_{n-1}$ o\`u $\sigma$ est 
d\'efini, on a \ $Q\in\sigma(P)$ \ si et seulement si \ 
$\phi(\sigma(P)).\sigma(Q) = 0$ .

2 - L'image de $\sigma$ est dense dans \m{Q(\phi)}, et
\m{\phi_S} est de rang au moins 3.
\rm
\end{subsub}

\begin{proof} La d\'emonstration de 1- est la m\^eme que celle de \cite{dr3}, 
lemme 2.5. Puisque l'image de \m{Q(\phi)} n'est pas contenue dans un hyperplan,
2- en d\'ecoule imm\'ediatement.
\end{proof}

\vskip 1.5cm

\begin{sub}\label{norm} \bf Congruences quadratiques normales\rm

Le cas le plus int\'eressant est celui o\`u $W_0$ est l'espace des formes
quadratiques s'annulant sur une sous-vari\'et\'e ferm\'ee $C$ de \m{\P_{n-1}}
de codimension 2 en chacun de ses points. Dans ce cas $C$ est n\'ecessairement 
une quadrique d'un hyperplan de \m{\P_{n-1}}, en vertu de la

\vskip 1cm

\begin{subsub}\label{intersquad}{\bf Proposition : } Soit $C$ une 
sous-vari\'et\'e ferm\'ee de \m{\P_{n-1}} de codimension 2 en chacun de ses 
points. Alors les 2 assertions suivantes sont \'equivalentes :
\begin{enumerate}
\item Il existe un hyperplan $H$ de $\P_{n-1}$ tel que $C$ soit une quadrique
de $H$.
\item L'id\'eal de $C$ contient $n+1$ formes quadratiques lin\'eairement
ind\'ependantes.
\end{enumerate}
\end{subsub}

\begin{proof}
Je pense que ce r\'esultat est classique. J'en donne cependant une esquisse de
d\'emonstration. Il est imm\'ediat que 1- entraine 2-. Prouvons la r\'eciproque.
Supposons que 2- soit vraie. 

Soit \ \m{W=H^0(\ki_Z(2))}, qui est de dimension
\m{\geq n+1}. Supposons que \m{\P(W)} contienne une quadrique $Q$ de rang
\m{r\geq 5}. Alors $C$ est une hypersurface de $Q$. D'apr\`es le th\'eor\`eme de
Klein (cf. \cite{ha}, Ex. II.6.5), $C$ est l'intersection compl\`ete de $Q$ et
d'une hypersurface $S$ de \m{\P_{n-1}}. L'id\'eal de $C$ est donc engendr\'e par
des \'equations $q$ et $\phi$ de $Q$, $S$ respectivement. Il en d\'ecoule 
d'apr\`es 2-
que $\phi$ doit \^etre de degr\'e 1, ce qui d\'emontre 1-. 

On peut donc supposer que toutes les formes quadratiques de $W$ sont de rang
\m{\leq 4}. Soit $r$ leur rang maximal. On voit ais\'ement que les quadriques
associ\'ees ont un \m{\P_{n-r}} en commun, et on se ram\`ene ainsi au cas o\`u \
\m{r=n} \ et \ \m{n=3} ou $4$. Traitons par exemple le cas le plus difficile :
\m{n=4}. On prend une quadrique lisse \ \m{Q\simeq\P_1\times\P_1} \ de 
\m{\P(W)}. Alors $C$ est une section d'un fibr\'e \m{\ko(a,b)} dans $Q$, avec \
\m{a+b\leq 4} (car $C$ est contenue dans l'intersection de deux quadriques). On
montre ensuite que le seul cas possible est \ \m{(a,b)=(1,1)}, en remarquant que
dans les autres cas l'id\'eal de $C$ ne contient pas 5 formes quadratiques
lin\'eairement ind\'ependantes. Ceci d\'emontre 1-.
\end{proof}

\bigskip

Les congruences quadratiques \ \m{\P_{n-1}\lra\P(W)} \ sont dites {\em
normales}. Ce sont les seuls types de congruences quadratiques que l'on
\'etudiera ici. Il en existe d'autres (cf. \ref{congplan}).

Soient $C$ une quadrique d'un hyperplan $H$ de \m{\P_{n-1}}, $W$ l'espace
vectoriel de dimension \m{n+1} des quadriques de \m{\P_{n-1}} contenant $C$.
Soient \ \m{\sigma : \P_{n-1}\lra\P(W)} \ une congruence quadratique, \m{\phi :
W\simeq W^*} \ l'isomorphisme associ\'e, et $r$ le rang de \m{Q(\phi)}. On 
notera \ \m{H(\sigma)=H} \ et \ \m{C(\sigma)=C}.

On pose \ \m{m=rg(C(\sigma))}. On peut supposer que \m{\P_{n-1}} est muni des
coordonn\'ees \m{x_1,\ldots,x_n}, que \m{H(\sigma)} est d\'efini par 
l'\'equation \ \m{x_n=0}, et \m{C(\sigma)} est d\'efini dans \m{H(\sigma)} par
l'\'equation \ \m{x_1^2+\cdots+x_m^2=0}. On munit $W$ de la base 
\m{(x_1x_n,\ldots,x_n^2,x_1^2+\cdots+x_m^2)}. On consid\`ere la matrice 
\m{(n+1)\times(n+1)}

\[
J_m \ = \ \left(\begin{array}{ccccccccccc}
1 & 0 & . & . & . & . & . & . & . & 0 & 0 \\
0 & 1 & 0 & . & . & . & . & . & . & . & 0 \\
. & . & . & . & . & . & . & . & . & . & 0 \\
. & . & . & . & . & . & . & . & . & . & 0 \\
0 & . & . & . & . & 1 & 0 & . & . & . & 0 \\
0 & . & . & . & . & 0 & 0 & . & . & . & 0 \\
. & . & . & . & . & . & . & . & . & . & . \\
. & . & . & . & . & . & . & . & . & . & . \\
0 & . & . & . & . & . & . & . & . & 0 & 0 \\
0 & . & . & . & . & . & . & . & 0 & 0 & -\frac{1}{2} \\ 
0 & . & . & . & . & . & . & . & 0 & -\frac{1}{2} & 0 \\ 
\end{array}\right)
\]
dont les seuls termes non nuls de la diagonale sont les $m$ premiers termes,
\'egaux \`a $1$. Alors on voit ais\'ement qu'il existe une matrice 
\m{(n+1)\times(n+1)} antisym\'etrique $M$ telle que la matrice de $\phi^{-1}$
relativement \`a la base pr\'ec\'edente de $W$ soit, \`a un scalaire
multiplicatif pr\`es de la forme \ \m{M+J_m}.
Le rang de \m{\phi_S^{-1}} est alors \m{m+2}, ainsi donc que celui de 
\m{\phi_S}. On en d\'eduit imm\'ediatement la

\vskip 1cm

\begin{subsub}\label{lemm5}{\bf Proposition : } Le rang de \m{C(\sigma)} est 
\m{r-2}, \m{(\phi^{-1})_S} ne d\'epend que de $W$ et $\sigma$ induit un
isomorphisme birationnel \ \m{\P_{n-1}\simeq\ov{\sigma(\P_{n-1})}}.
\end{subsub}

\bigskip

On d\'eduit aussi de ce qui pr\'ec\`ede le

\vskip 1cm

\begin{subsub}\label{lemm6}{\bf Lemme : } La restriction de $\sigma$ \`a
\m{H(\sigma)} est constante, et \m{\sigma(H(\sigma))} est une quadrique 
d\'eg\'en\'er\'ee, dont une des composantes est \m{H(\sigma)}.
\end{subsub}

\medskip

On notera \m{L(\sigma)} l'autre composante de \m{\sigma(H(\sigma))}. C'est un
hyperplan de \m{\P_{n-1}}. Avec les notations pr\'ec\'edentes, si
\ \m{M=(a_{ij})}, l'\'equation de \m{L(\sigma)} est
\[\sigg_{i=1}^n a_{n+1,i}X_i \ - \ \frac{1}{2}X_n \ \ = \ \ 0.\]

\end{sub}

\vskip 1.5cm

\begin{sub}\label{carac}\bf Caract\'erisation des congruences quadratiques\rm

\begin{subsub}\label{propcarac}{\bf Proposition : } Soient $W_0$ un sous-espace
vectoriel de dimension \m{n+1} de 

\noindent\m{H^0(\P_{n-1},\ko(2))} et \ \m{\sigma:\P_{n-1}\to\P(W_0)} \ un 
morphisme rationnel dont l'image n'est pas contenue dans un hyperplan, et tel 
que pour tout \m{P\in\P_{n-1}} g\'en\'eral on ait \m{P\in\sigma(P)}. Si 
$P\in\P_{n-1}$, on note $F_\sigma(P)$ l'adh\'erence du lieu des points $Q$ tels
que $\sigma$ soit d\'efini en $Q$ et $P\in\sigma(Q)$. Alors $\sigma$ est une 
congruence quadratique si et seulement si pour tout point g\'en\'eral $P$ de 
$\P_{n-1}$, $F_\sigma(P)$ est un point de de $\P(W_0)$.
\end{subsub}

\begin{proof}
Posons \ $L=\sigma^*(\ko(1))$, et soit \ \m{\sigma^*:W_0^*\lra H^0(L)} \ 
l'application lin\'eaire associ\'ee \`a $\sigma$. 
Soient $x\in V$ au dessus de $P$, et $\phi_x\in W_0^*$ d\'efini par 
\ \m{\phi_x(q)=q(x)}. Alors $\sigma^*(\phi_x)$ s'annule exactement sur
$F_\sigma(P)$. Il en d\'ecoule d'apr\`es la condition 2 de \ref{cong}
que si $\sigma$ est une congruence quadratique,
alors $F_\sigma(P)$ est une quadrique \'el\'ement de $\P(W_0)$.

R\'eciproquement, supposons que pour un point g\'en\'eral $P$ de $\P_{n-1}$, 
$F_\sigma(P)$ soit une quadrique appartenant \`a $\P(W_0)$. Il faut prouver que
les propri\'et\'es 1- et 2- de \ref{cong} sont v\'erifi\'ees. En choisissant 
\m{n+1} points \m{x_1,\ldots, x_{n+1}} de $V$ tels que \m{(\phi_{x_i})} soit 
une base de $W_0^*$ on se ram\`ene au cas o\`u relativement
\`a une base convenable de $W_0$, $\sigma$ est de la forme
\[
\sigma \ = \ (\phi_1^p,\ldots,\phi_{n+1}^p)
\]
$p$ \'etant un entier positif et \m{\phi_1,\ldots,\phi_{n+1}\in W_0}. 

Montrons que \m{p=1}. Le morphisme 
\[
\begin{array}{ccc}
\P(V) & \lra & \P(W_0^*) \\ {\C}x & \longmapsto & {\C\phi_x}
\end{array}
\]
a pour image une quadrique. Il en d\'ecoule qu'il existe une quadrique Q de
\m{\P_n} telle que pour tout \ \m{(\lambda_1,\ldots,\lambda_{n+1})
\in Q} \ il existe un \m{\phi\in W_0} tel que
\[
\lambda_1\phi_1^p+\cdots+\lambda_{n+1}\phi_{n+1}^p \ =
\ \phi^p.
\]
Puisque l'image de $\sigma$ n'est pas contenue dans un hyperplan, \m{\phi_1^p,
\ldots, \phi_{n+1}^p} sont lin\'eairement ind\'ependants
dans \m{S^{2p}V^*}. La quadrique $Q$ contient des droites. Cela signifie qu'il
existe des \'el\'ements lin\'eairements ind\'ependants $\psi_1$, $\psi_2$ de
$W_0$ tels que toute combinaison lin\'eaire de $\psi_1^p$ et $\psi_2^p$ soit de
la forme $\psi^p$, avec $\psi$ dans $W_0$. Ceci est impossible si $p>1$, comme
on peut le voir par exemple en se restreignant \`a des droites de $\P_{n-1}$.
Il est maintenant imm\'ediat que les conditions 1- et 2- de \ref{cong} sont
v\'erifi\'ees. 
\end{proof}

\vskip 1cm

\begin{subsub}\label{propcarac2}{\bf Proposition : } Avec les notations de la
proposition \ref{propcarac}, le morphisme rationnel
\[
\begin{array}{ccc}
\P_{n-1} & \lra & \P(W) \\ P & \longmapsto & F_\sigma(P)
\end{array}
\]
est une congruence quadratique.
\end{subsub}

\begin{proof}
Notons $\tau$ le morphisme pr\'ec\'edent. Il suffit de v\'erifier que les
conditions de la proposition \ref{propcarac} sont v\'erifi\'ees par $\tau$.
Il est imm\'ediat que pour un point g\'en\'eral $P$ de \m{\P_{n-1}} on a
\ \m{F_\tau(P)=\sigma(P)}, donc \m{F_\tau(P)} est bien une quadrique. 
D'autre part,
si \m{x\in V} est au dessus de $P$ on a vu que \m{\sigma^*(\phi_x)} est une
\'equation de  \m{F_\sigma(P)}. Comme \m{\sigma^*} est un isomorphisme de
\m{W_0^*} sur \m{W_0} et que les \m{\phi_x} engendrent \m{W_0^*}, on voit que
l'image de $\tau$ ne peut pas \^etre contenue dans un hyperplan de \m{\P(W_0)}.
\end{proof}

\vskip 1cm

\begin{subsub}\label{deftrans}\rm
La congruence quadratique de la proposition \ref{propcarac2} est appel\'ee la
{\em transpos\'ee} de $\sigma$ et est not\'ee \ \m{^t\sigma}. Notons que
l'isomorphisme \ \m{W_0\simeq W_0^*} \ associ\'e \`a $^t\sigma$ est le 
transpos\'e
de celui qui est associ\'e \`a $\sigma$. On a \ \m{^t(^t\sigma)=\sigma}.
\end{subsub} 

\end{sub}

\vskip 1.5cm

\begin{sub}\label{transl}\bf Translation associ\'ee\rm

Soit $\sigma$ une congruence quadratique induisant un isomorphisme birationnel
\ \m{\P_{n-1}\simeq\ov{\sigma(\P_{n-1})}}.
Soit \  \m{\phi=(\sigma^*)^{-1}}. On note
\[
\T(\sigma) \ = \ \sigma^{-1}\circ\T(\phi)\circ\sigma,
\]
qu'on appelle la {\em translation} associ\'ee \`a $\sigma$.  

\vskip 1cm

\begin{subsub}\label{lemm2}{\bf Lemme : } Soit $P$ est un point de \m{\P_{n-1}}  
o\`u
$\sigma$ est d\'efinie. Alors $\sigma$ est aussi d\'efinie en \m{\T(\sigma)(P)}
et \ \m{\T(\sigma)(P)\in\sigma(P)}. \end{subsub}

\begin{proof}
Imm\'ediat. 
\end{proof}

\vskip 1cm

\begin{subsub}\label{lemm2b}{\bf Proposition : } 
Soient $\sigma$ une congruence quadratique r\'eguli\`ere
ou normale, \hfil\break\m{\phi : 
W_0\simeq W_0^*} l'isomorphisme associ\'e. Alors \m{^t\sigma} est
aussi r\'eguli\`ere, et on a
\[ \T(^t\sigma) \ = \ \T(\sigma)^{-1}. \]\end{subsub}

\begin{proof} Le cas o\`u $\sigma$ est normale d\'ecoule du \para 3.5. Supposons
donc $\sigma$ r\'eguli\`ere.
Soit \ \m{\phi = (\sigma^*)^{-1} : W_0\simeq W_0^*}. Si $P$ est un point de
\m{\P_{n-1}} o\`u $\sigma$ est d\'efinie, on note \m{W(P)} le sous-espace 
vectoriel de \m{W_0} constitu\'e de $u$ tels que \ \m{\phi(u)\sigma(P)=0}. 
Puisque \m{\phi_S} est de rang au moins $n$ et l'image de $\sigma$ dense dans 
\m{Q(\phi)} (d'apr\`es le lemme \ref{lemm1}), 
pour un $P$ g\'en\'erique \m{W(P)} est un hyperplan de \m{W_0} et \m{\phi_S} est
non d\'eg\'en\'er\'ee sur \m{W(P)}.
Il faut montrer que \m{^t\sigma} est injective sur un ouvert de \m{\P_{n-1}}.
Supposons que ce ne soit pas le cas. Alors il existe deux points distincts
\m{Q_0},
\m{Q_1} de \m{\P_{n-1}} o\`u $\sigma$ est d\'efinie, tels que pour \ \m{i=0,1}
\m{W(Q_i)} soit un hyperplan de \m{W_0}, que \m{\phi_S} soit non
d\'eg\'en\'er\'ee sur 
\m{W(Q_i)} et \ \m{{^t\sigma}(Q_0)= {^t\sigma}(Q_1)}. Cela entraine que les
quadriques de la forme \m{\sigma(Q)} qui passent par \m{Q_0} sont les m\^emes
que celles qui passent par \m{Q_1}. D'apr\`es le lemme \ref{lemm1}, ceci
\'equivaut \`a \ \m{Q(\phi)\cap W(Q_0) = Q(\phi)\cap W(Q_1)}. Puisque $\phi_S$
est non d\'eg\'en\'er\'ee sur \m{W(Q_0)} et \m{W(Q_1)}, on en d\'eduit que \
\m{W(Q_0)=W(Q_1)}. Puisque \m{\phi_S} est non d\'eg\'en\'er\'ee, ceci entraine 
que \ \m{Q_0=Q_1}. Donc \m{^t\sigma} est injective sur un ouvert de 
\m{\P_{n-1}}. Le reste de la proposition \ref{lemm2b} est \'evident.
\end{proof}
\end{sub}

\vskip 1.5cm

\begin{sub}\label{constr}{\bf Construction g\'eom\'etrique des congruences 
quadratiques normales}\rm

Soient $n$ un entier tel que \ \m{n\geq 3}, $W$ un $\C$-espace vectoriel de 
dimension \ \m{n+1}, \ \m{\phi : W\lra W^*} \ un isomorphisme. On suppose que
la quadrique \m{Q(\phi)} est non d\'eg\'en\'er\'ee. Soit \ \m{r\geq 3} \ son
rang. Si \ \m{X\subset\P(W)}, rappelons qu'on note \m{X^\bot} (resp. 
\m{^\bot X}) le sous-espace lin\'eaire de \m{\P(W)} constitu\'e des points 
\m{\C y} tels que \ \m{x\phi(y)=0} (resp. \ \m{y\phi(x)=0}) pour tout point $x$
de $W$ au dessus d'un point de $X$.

Soient $O$ un point lisse de \m{Q(\phi)} et 
\[ \pi = \pi_O : \P(W)\backslash\lbrace O\rbrace\lra\P_{n-1} \]
la projection de centre $O$. Plus concr\`etement, on peut consid\'erer que
\m{\P_{n-1}} est un hyperplan de \m{\P(W)} ne contenant pas $O$ et pour tout
point $P$ de \m{\P(W)} diff\'erent de $O$, \m{\pi(O)} est l'intersection de
\m{\P_{n-1}} et de la droite \m{OP}. 

\vskip 1cm

\begin{subsub}\label{inters1}{\bf Lemme : } Soit \m{T_O} l'hyperplan tangent 
\`a \m{Q(\phi)} en $O$. L'image par $\pi$ de \ \m{Q(\phi)\cap T_O} est une 
quadrique $C$ de l' hyperplan \m{\pi(T_O)} de \m{\P_{n-1}}. Le rang de $C$ 
est \'egal \`a \m{r-2}.
\end{subsub}

\begin{proof}
On peut choisir des coordonn\'ees ind\'ependantes \m{x_0,\ldots,x_n} dans
\m{\P(W)} de telle sorte que l'\'equation de \m{Q(\phi)} soit \
\m{x_0^2+\cdots+x_{r-1}^2=0}, et que $O$ soit le point \m{(1,i,0,\ldots,0)}.
L'\'equation de \m{T_O} est alors \ \m{x_0+ix_1=0}. On peut aussi supposer que
\m{\P_{n-1}} est l'hyperplan d'\'equation \ \m{x_0=0}, muni des coordonn\'ees
\m{x_1,\ldots,x_n}. Dans ce cas, \m{\pi(T_O)} est l'hyperplan de \m{\P_{n-1}}
d'\'equation \ \m{x_1=0}, et $C$ la quadrique d'\'equation \
\m{x_2^2+\cdots+x_{r-1}^2=0}, qui est bien de rang \m{r-2}. 
\end{proof}

\bigskip

On note \m{W_0} l'espace vectoriel de dimension \m{n+1} des quadriques de
\m{\P_{n-1}} contenant $C$. Il est imm\'ediat que si $H$ est un hyperplan de 
\m{\P(W)}, \m{\pi(Q(\phi)\cap H)} \ est une quadrique de \m{\P_{n-1}} contenant 
$C$, qu'on notera \m{S(H)}. On a des isomorphismes
\[
\begin{array}{cccc}
\tau : & \P(W) & \lra & \P(W_0) \\
       & P     & \longmapsto & S(P^\bot)
      
\end{array}
\]
et
\[
\begin{array}{cccc}
\tau' : & \P(W) & \lra & \P(W_0) \\
       & P     & \longmapsto & S({^\bot}P)
      
\end{array}
\]

\medskip

{\bf Notations : }\rm On notera pour simplifier de la m\^eme fa\c con un
point $P$ d'un espace projectif, et un point de l'espace vectoriel correspondant
au dessus de $P$. Si $u$ est un \'el\'ement d'un espace vectoriel et $f$ un
\'el\'ement de son dual, on notera \m{u.f} ou \m{f.u} le scalaire image de $u$
par $f$.

\medskip

Si $Q$ est un point de \ \m{\P_{n-1}\backslash C}, on note \m{\eta(Q)} le point
commun \`a la droite \m{OQ} et \`a \m{Q(\phi)} autre que $O$. On a
\begin{eqnarray*}
\eta(Q) & = & (Q.\phi(Q))O-(Q.\phi(O)+O.\phi(Q))Q \\
        & = & (Q.\phi_S(Q))O-(Q.\phi_S(O)+O.\phi_S(Q))Q. \\
\end{eqnarray*}
Si \m{P\in\P(W)}, on a donc \m{Q\in\tau(P)} si et seulement si \ 
\m{P.\phi(\eta(Q))=0}, et \m{Q\in\tau'(P)} si et seulement si \ 
\m{P.{^t\phi}(\eta(Q))=0}.

\medskip

On consid\`ere les morphismes rationnels
\[
\begin{array}{cccc}
\sigma : & \P_{n-1} & \lra        & \P(W_0)        \\
         & Q        & \longmapsto & S({^\bot\eta(Q)})
\end{array}
\]
et
\[
\begin{array}{cccc}
\sigma' : & \P_{n-1} & \lra        & \P(W_0)          \\
          & Q        & \longmapsto & S(\eta(Q)^\bot)
\end{array}
\]
Si $P$, $Q$ sont des points g\'en\'eraux de \m{\P_{n-1}} on a donc \
\m{P\in\sigma(Q)} (resp. \ \m{P\in\sigma'(Q)}) si et seulement si \
\m{\eta(P).\phi(\eta(Q))=0} (resp. \m{\eta(Q).\phi(\eta(P))=0}).

\vskip 1cm

\begin{subsub}\label{prop5}{\bf Proposition : } Les morphismes $\sigma$ et
\m{\sigma'} sont des congruences quadratiques. On a \ \m{\sigma'={^t \sigma}} \
et l'isomorphisme \ \m{W_0\simeq W_0^*} \ associ\'e \`a $\sigma$ est $\phi$.
\end{subsub}

\begin{proof} Pour tout point $P$ de \m{\P_{n-1}}, on note \m{\Phi_P} la forme
lin\'eaire \ \m{q\longmapsto q(P)}. On commence par \'evaluer 
\m{\sigma^*(\Phi_P)}, pour tout point $P$ de \m{\P_{n-1}\backslash C}. Soit
\m{Q\in\P_{n-1}}. Alors on a
\[\sigma^*(\Phi_P)(Q) \ = \ \eta(P).\phi(\eta(Q)).\]
Donc \ \m{\tau(\eta(P))=\sigma^*(\Phi_P)}. On a d'autre part, pour tout \
\m{w\in W}, \ \m{^t\tau(\Phi_P)(w)=w.\phi(\eta(P))}. Il en d\'ecoule, les
\m{\Phi_P} engendrant \m{W_0^*}, que \m{\sigma^*} est \`a valeurs dans \m{W_0}
et que \ \m{^t\tau\circ(\sigma^*)^{-1}\circ\tau=\phi}. La proposition 
\ref{prop5} s'en d\'eduit imm\'ediatement.
\end{proof}

\bigskip

On a donc \ \m{H(\sigma)=\pi(T_O)} \ et \ \m{C(\sigma)=C}.

\vskip 1cm

\begin{subsub}\label{lemm7}{\bf Proposition : } On a \ \m{L(\sigma)=\pi(O^\bot)}.
\end{subsub}

\begin{proof} Si $P$ est un point g\'en\'eral de \m{H(\sigma)}, on a \
\m{\eta(P)=O}, donc \ \m{\pi(O^\bot)\subset\sigma(P)}. Si \ \m{T_0\not= O^\bot}
\ on a donc \ \m{L(\sigma)=\pi(O^\bot)}. Le cas \ \m{T_0= O^\bot} \ s'en
d\'eduit par continuit\'e.
\end{proof}

\bigskip

La translation \m{\T(\sigma)} s'interpr\`ete de la fa\c con suivante : si
\m{P\in\P_{n-1}}, on a
\[\T(\sigma)(P) \ = \ \pi(\T(\phi)(\eta(P))). \]
On en d\'eduit la

\vskip 1cm

\begin{subsub}\label{prop6}{\bf Proposition : } Les propri\'et\'es suivantes
sont \'equivalentes :
\begin{enumerate}
\item On a \ $T_0=O^\bot$.
\item On a \ $L(\sigma)=H(\sigma)$.
\item La translation $T(\sigma)$ est lin\'eaire.
\end{enumerate}
\end{subsub}

\begin{proof} L'\'equivalence des deux premi\`eres propri\'et\'es
pr\'ec\'edentes d\'ecoule du lemme \ref{lemm7}. Supposons que 1- soit
v\'erifi\'ee. Soit $K$ un hyperplan de \m{\P_{n-1}}. Alors les points
\m{\eta(Q)}, \m{Q\in K} sont dans \m{K'\cap Q(\phi)}, o\`u \ \m{K'=\pi^{-1}(K)}
\ est un hyperplan de \m{\P(W)} passant par $O$. On a \ \m{\T(\sigma)(K) \ = \ 
\pi(\T(\phi)(K'\cap Q(\phi))}, et \m{\T(\phi)(K'\cap Q(\phi))} est contenu dans
\m{\T(\phi)(K')}, qui est un hyperplan de \m{\P(W)}. Son image par $\pi$ est 
contenue dans un hyperplan de \m{\P_{n-1}} car puisque \ \m{T_0=O^\bot}, on a \
d'apr\`es le lemme \ref{lemm8}, \m{\T(\phi)(O)=O}, d'o\`u \m{O\in\T(\phi)(K')}.
Il en d\'ecoule que \m{\T(\sigma)} est lin\'eaire. R\'eciproquement, si   
\m{\T(\sigma)} est lin\'eaire, le raisonnement pr\'ec\'edent montre que pour
tout hyperplan $K'$ de \m{\P(W)} passant par $O$ g\'en\'eral, \m{\T(\phi)(K')}
passe aussi par $O$, d'\`ou il d\'ecoule que \ \m{\T(\phi)(O)=O} \ dans
\m{\P(W)} et \ \m{T_0=O^\bot} \ d'apr\`es le lemme \ref{lemm8}.
\end{proof} 

\vskip 1cm

\begin{subsub}\label{theo3}{\bf Th\`eor\`eme : } Toute congruence quadratique
normale peut s'obtenir par la m\'ethode pr\'ec\'edente.
\end{subsub}

\begin{proof}
On emploie des notations l\'eg\`erement diff\'erentes de celles de \ref{norm}.
On munit \m{\P_n} de coordonn\'ees ind\'ependantes \m{x_1,\ldots,x_{n+1}}. On
suppose que \m{\P_{n-1}\subset\P_n} est l'hyperplan d'\'equation \m{x_{n+1}=0},
muni des coordonn\'ees \m{x_1,\ldots,x_n}, 
que \m{Q(\phi)} est la quadrique d'\'equation
\[ x_{n-m+1}^2+\cdots+x_{n+1}^2 \ = \ 0, \]
et que \ \m{O=(0,\ldots,i,1)}. Alors \m{T_O} est l'hyperplan d'\'equation \
\m{ix_n+x_{n+1}=0}. L'hyperplan \m{\pi(T_O)} de \m{\P_{n-1}} est d\'efini par
l'\'equation \ \m{x_n=0} \  et $C$ est la quadrique de \m{\pi(T_O)} d'\'equation
\[ x_{n-m+1}^2+\cdots+x_{n-1}^2 \ = \ 0. \]
L'isomorphisme $\phi$ a une matrice du type
\[
\left(\begin{array}{cc} 0 & 0 \\ 0 & I_{m+1}\\
\end{array}\right) \ + \ A,
\]
o\`u \m{I_{m+1}} est la matrice identit\'e \m{(m+1)\times(m+1)} et 
\ \m{A=(\alpha_{ij})} \ une
matrice antisym\'etrique. Si \ \m{P=(x_1,\ldots,x_n)\in\P_n}, on calcule
ais\'ement qu'on a, en posant \ \m{s(P)=x_{n-m+1}^2+\cdots+x_{n-1}^2},
\[
\eta(P) \ = \ \left(\begin{array}{c} -2ix_1x_n \\ . \\ . \\ . \\
-2ix_{n-1}x_n \\ is(P)-ix_n^2 \\ s(P)+x_n^2 \\ \end{array}\right) .
\]
On va montrer que toute congruence quadratique \ \m{\sigma :
\P_{n-1}\lra\P(W_0)} \ peut s'obtenir par la construction g\'eom\'etrique
d\'ecrite pr\'ec\'edemment en utilisant un $\phi$ convenable (c'est-\break
\`a-dire une
matrice antisym\'etrique $A$ convenable). On munit \m{W_0} de la base
\[
(p_1,\ldots,p_{n+1}) \ = \ 
(x_1x_n,\ldots,x_{n-1}x_n,x_n^2,x_{n-m+1}^2+\cdots+x_{n-1}^2).
\]
Soit 
\[N \ = \ \left(\begin{array}{cc} 0 & 2 \\ 2 & 0 \\
\end{array}\right)
\]
Alors $\sigma$ est d\'efinie par une matrice \m{(m+1)\times(m+1)} 
 \ \m{B=(b_{ij})} de la forme
\[
B \ = \ \left(\begin{array}{cc}-4I_{m-1} & 0 \\ 0 & N \\ \end{array}\right),
\]
c'est-\`a-dire que si $P$ et $Q$ sont des 
points de \m{\P_{n-1}}, on a \ \m{P\in\sigma(Q)} si et seulement si
\[\sigg_{1\leq k,j\leq n+1}b_{jk}p_k(P)p_j(Q) \ = \ 0.\]
En utilisant le fait que \m{P\in\sigma(Q)} si et seulement si \
\m{\eta(Q).\phi(\eta(P))=0}, on voit qu'on doit avoir
\[ b_{ij}=-4\alpha_{kj} {\rm \ \ \ \ si \ \ \ \ } 1\leq k,j\leq n-1, \]
\[ b_{kn} = -2\alpha_{kn}-2i\alpha_{k,n+1} {\rm \ \ \ \ si \ \ \ \ }
1\leq k\leq n-1, \]
\[ b_{k,n+1} = 2\alpha_{kn}-2i\alpha_{k,n+1}  {\rm \ \ \ \ si \ \ \ \ }
1\leq k\leq n-1, \]
\[ b_{n,n+1} = 2i\alpha_{n,n+1}. \]
Il est donc clair qu'on peut bien choisir une matrice antisym\'etrique $A$
ad\'equate. 
\end{proof}
\end{sub}

\vskip 1.5cm

\begin{sub}{\bf \'Equivalence de congruences quadratiques}\rm

Soit \ \m{\sigma :\P_{n-1}=\P(V)\lra\P_n=\P(W_0)} \ une congruence quadratique. 
Donc \m{W_0} est un sous-espace vectoriel de dimension \m{n+1} de 
\m{H^0(\P_{n-1},\ko(2))}. On note \m{G_{W_0}} le sous-groupe de \m{GL(V)}
constitu\'e des \'el\'ements qui laissent \m{W_0} invariant. On note 
\m{\Cong(W_0)} l'ensemble des congruences quadratiques \`a valeurs dans 
\m{\P(W_0)}. On a une action de \m{G_{W_0}} sur \m{\Cong(W_0)} d\'efinie par :
\[ (\alpha.\sigma)(P) \ = \ \alpha^{-1}(\sigma(\alpha(P))) \]
pour \m{P\in\P_{n-1}} et \m{\alpha\in G_{W_0}}. 

Le morphisme
\begin{eqnarray*}
\Phi_{W_0} : \Cong(W_0) & \lra & \Isom(W_0,W_0^*)\\
\sigma & \longmapsto & (\sigma^*)^{-1}=\phi\\
\end{eqnarray*}
est compatible avec le morphisme de groupes \'evident \ \m{G_{W_0}\lra GL(W_0)}.

Si \m{W_0} est l'espace des \'equations de quadriques contenant une quadrique
$C$ d'un hyperplan de \m{\P_{n-1}} (autrement dit si on s'int\'eresse \`a des
congruences quadratiques normales), \m{G_{W_0}} est le groupe des \'el\'ements
de \m{GL(V)} laissant $C$ invariante. On notera alors
\[ G_C=G_{W_0}, \ \ \ \Cong(C)=\Cong(W_0), \ \ \ \Phi_C=\Phi_{W_0}.\] 
\end{sub}

On montre au \para \ref{equiv} que si \ \m{n=3}, les fibres du morphisme
quotient
\[\Cong(C)/G_C\lra\Isom(W_0,W_0^*)/GL(W_0)\]
sont finies. J'ignore si c'est le cas si \ \m{n>3}. Il est vraisemblable que non
car si \m{n>3}, et \ \m{\phi\in\Isom(W_0,W_0^*)}, on a en g\'en\'eral \
\m{\dim(G_\phi)<\dim(Q(\phi))} (cf. \cite{go3}, 3.5.2), \m{G_\phi} d\'esignant
le stabilisateur de $\phi$ dans \m{GL(W_0)}, c'est-\`a-dire le groupe des {\em
isom\'etries} de la forme bilin\'eaire induite par $\phi$. 

\vskip 2.5cm

\section{Congruences quadratiques planes}\label{congplan}

On reprend les notations de \ref{congdef}. 

\vskip 1cm

\begin{sub}\label{def2}\bf Classification des congruences quadratiques\rm

\begin{subsub}{\bf Proposition : }\label{prop1}
Soit \ \m{\sigma:\P_2\lra\P(W)} \ une congruence quadratique. Alors
on est dans l'un des trois cas suivants :

\noindent 1) $\P(W)$ est l'espace des coniques passant par deux points
fixes distincts de $\P_2$.

\noindent 2) $\P(W)$ est l'espace des coniques passant par un point fixe 
de $\P_2$ et tangentes \`a une droite fixe de $\P_2$ en ce point.

\noindent 3) $\P(W)$ est l'espace des coniques invariantes par une 
involution non triviale de $\P_2$.

Dans les cas 1 et 2, $\sigma$ induit un isomorphisme birationnel 
 \ $\P_2\simeq Q(\phi)$. \rm

(cf. \cite{dr3}, prop. 2.6). 

\bigskip

Dans le cas 1 (resp. 2,3) on dit que $\sigma$ est
une {\em congruence quadratique de type} 1 (resp. 2,3).
Dans les cas 1 et 2 on est en pr\'esence de congruences quadratiques {\em
normales} (cf. \ref{norm}). Il est facile de trouver des exemples du cas 3, qui
montrent qu'il existe des congruences quadratiques non normales.

Dans ce qui suit on s'int\'eresse exclusivement aux cas 1 et 2.
\end{subsub}

Soit
\[
\pi : S^2W\lra H^0(\P_2,\ko(4))
\]
la restriction du morphisme canonique \ \m{S^2(H^0(\P_2,\ko(2)))\lra 
H^0(\P_2,\ko(4))}. Dans les 3 cas de la proposition pr\'ec\'edente,
\m{\ker(\phi)} est de dimension 1. Un g\'en\'erateur $\psi$ de \m{\ker(\phi)} 
est une \'equation de l'image de \m{Q(\phi)} par l'isomorphisme \
\m{\P(W)\simeq\P(W^*)} \ d\'eduit de $\phi$ (c'est la quadrique de \m{\P(W^*)}
d'\'equation \ \m{x.\phi^{-1}(x)=0}). On a aussi
\[\psi \ = \ (\phi^{-1})_S.
\]
On en d\'eduit le r\'esultat suivant, qui d\'ecoule aussi du th\'eor\`eme
\ref{theo3} :

\vskip 1cm

\begin{subsub}{\bf Proposition : }\label{prop2} Soit \ \m{\phi\in\Isom(W,W^*)}. 
Alors

\noindent 1 - $\phi$ admet une r\'ealisation g\'eom\'etrique par des coniques de
type 1 si et seulement si \m{Q(\phi)} est de rang 4.

\noindent 2 - $\phi$ admet une r\'ealisation g\'eom\'etrique par des coniques de
type 2 si et seulement si \m{Q(\phi)} est de rang 3.
\end{subsub}

\begin{proof}
Soit \ \m{\sigma : \P_2\lra\P(W)} \ un congruence quadratique de type 1 ou 2.
On peut donc trouver 
des coordonn\'ees ind\'ependantes $X$, $Y$, $Z$ dans $V$ de telle sorte
qu'on ait la base suivante \m{(e_0,e_1,e_2,e_3)} dans $W$ :

\noindent Cas 1 : \m{e_0=XY}, \m{e_1=YZ}, \m{e_2=XZ}, \m{e_3=Z^2}, \
\m{\ker(\pi)=\pline{e_0e_3-e_1e_2}}.

\noindent Cas 2 : \m{e_0=XZ}, \m{e_1=Y^2}, \m{e_2=YZ}, \m{e_3=Z^2}, \
\m{\ker(\pi)=\pline{e_2^2-e_1e_3}}.

\noindent On en d\'eduit imm\'ediatement que si $\phi$ admet une r\'ealisation 
g\'eom\'etrique par des coniques de type 2, alors \m{Q(\phi)} est de rang 3
et si $\phi$ admet une r\'ealisation g\'eom\'etrique par des coniques de type 1,
alors \m{Q(\phi)} est de rang 4. 

R\'eciproquement, on munit $V$ de la base duale de \m{(X,Y,Z)}. On consid\`ere
les matrices
\[
A_1=\left(\begin{array}{cccc}
0 & 0 & 0 & 1 \\ 0 & 0 & -1 & 0 \\ 0 & -1 & 0 & 0 \\ 1 & 0 & 0 & 0
\end{array}\right), \ \ \ 
A_2=\left(\begin{array}{cccc}
0 & 0 & 0 & 0 \\ 0 & 0 & 0 & -1 \\ 0 & 0 & 2 & 0 \\ 0 & -1 & 0 & 0
\end{array}\right), \ \ \ 
\] 
Soit \ \m{\phi\in\Isom(W,W^*)}. Si \ \m{rg(Q(\phi))=4}, relativement \`a une 
base
convenable de $V$, la matrice de $\phi^{-1}$ s'exprime sous la forme de la somme
de $A_1$ et d'une matrice antisym\'etrique :
\[
\phi^{-1} \ = \ \left(\begin{array}{cccc}
0 & a & b & c+1 \\ -a & 0 & d-1 & e \\ -b & -d-1 & 0 & f \\ -c+1 & -e & -f & 0 
\end{array}\right). 
\]
On d\'efinit alors $\sigma$ par
\[
\begin{array}{ccc}
\sigma(\pline{x,y,z}) & = & \pline{(-ayz-bxz+(1-c)z^2)XY+(axy-(d+1)xz-ez^2)YZ\\
  &  & +(bxy+(d-1)yz-fz^2)XZ+((c+1)xy+eyz+fxz)Z^2}.
  \end{array}
\]
et on obtient ainsi une r\'ealisation g\'eom\'etrique de $\phi$ par des coniques
de type 1. L'autre cas est analogue, en utilisant la matrice \m{A_2}.
\end{proof}
\end{sub}

\newpage

\begin{sub}\label{transl2}\bf Translation associ\'ee\rm

Supposons que \ \m{rg(Q(\phi))=4}, et que la matrice de $\phi^{-1}$ relativement
\`a \m{XY}, \m{YZ}, \m{XZ}, \m{Z^2} soit 
\[
\phi^{-1} \ = \ \left(\begin{array}{cccc}
0 & a & b & c+1 \\ -a & 0 & d-1 & e \\ -b & -d-1 & 0 & f \\ -c+1 & -e & -f & 0 
\end{array}\right). 
\]
Soit \ \m{\theta=af+cd-be}. On a alors
\[
\T(\sigma)(x,y,z) \ = \ (L_1L_4, L_2L_3, L_2L_4),
\]
avec
\[
L_1 \ = \ (\theta+c+d+1)x+2ez, \ \ \ L_2 \ = \ -2bx+(\theta-c-d+1)z,
\]
\[
L_3 \ = \ (\theta-c+d-1)y-2fz, \ \ \ L_4 \ = \ 2ay+(\theta-d+c-1)z.
\]

\medskip

Supposons que \ \m{rg(Q(\phi))=3}, et que la matrice de $\phi^{-1}$ relativement
\`a \m{XZ}, \m{Y^2}, \m{YZ}, \m{Z^2} soit 
\[
\phi^{-1} \ = \ \left(\begin{array}{cccc}
0 & a & b & c \\ -a & 0 & d & e-1 \\ -b & -d & 2 & f \\ -c & -e-1 & -f & 0 
\end{array}\right). 
\]
Soient \ \m{\theta=af+cd-be}, \ \m{\Delta=\det(\phi^{-1})}. On a alors
\[
\T(\sigma)(x,y,z) \ = \ (V_1, L_1L_2, L_2^2),
\]
avec
\[
L_1 \ = \ (\theta+b)y+2cz, \ \ \ L_2 \ = \ -2ay+(\theta-b)z,
\]
\[
V_1 \ = \Delta xz -2(\theta d-2ae+2a+bd)y^2-4(\theta e-af+cd+b)yz
-2(\theta f-2ec-bf+2c)z^2.
\]

On en d\'eduit apr\`es quelques calculs le cas particulier suivant de la
proposition \ref{prop6} :

\vskip 1cm

\begin{subsub}\label{lemm3}{\bf Lemme : } La translation $\T(\sigma)$ est 
lin\'eaire si et seulement si on est dans un des deux cas suivants :
\begin{enumerate}
\item On a \ $rg(Q(\phi))=4$, $a=b=0$ et
\[ T(x,y,z) \ = \ ((d+1)((c+1)(d+1)x+2ez), (d-1)((c+1)(d-1)y-2fx),
(c-1)(d^2-1)z).\]
\item On a \ $rg(Q(\phi))=3$, $a=d=0$ et
\[ T(x,y,z) \ = \ (ex-4(e^2-1)y-2(f(e+1)+\frac{2ec}{b})z,
(e-1)y-\frac{2c}{b}z, (e+1)z).\]
\end{enumerate}
\end{subsub}
\end{sub}

\newpage

\begin{sub}\label{degen}\bf Lieu des coniques d\'eg\'en\'er\'ees\rm

Soit \ \m{\sigma : \P(V)\lra\P(W)} \ 
une congruence quadratique de type 1 ou 2. L'espace de coniques
$W$ contient des coniques d\'eg\'en\'er\'ees. On munit $W$ de la base d\'ecrite
dans la proposition \ref{prop2}. Pour le type 1, une conique d'\'equation \
\m{uXY+vYZ+wXZ+tZ^2=0} \ est d\'eg\'en\'er\'ee si et seulement si on a \
\m{u(vw-ut)=0}. Dans ce cas le lieu des coniques d\'eg\'en\'er\'ees est donc la
r\'eunion d'un plan et d'une quadrique de $\P(W)$. Pour le type 2, une conique 
d'\'equation \ \m{uXZ+vY^2+wYZ+tZ^2=0} \ est d\'eg\'en\'er\'ee si et seulement 
si on a \ \m{uv=0}. Dans ce cas le lieu des coniques d\'eg\'en\'er\'ees est 
la r\'eunion de deux plans de $\P(W)$.

On note $D(\sigma)$ le lieu des points $P$ de $P(V)$ tels que $\sigma(P)$ soit
d\'eg\'en\'er\'ee. 

Supposons que $\sigma$ soit de type 1. Alors en g\'en\'eral $D(\sigma)$ est la 
r\'eunion de 6 droites de $\P(V)$. On choisit des
coordonn\'ees ind\'ependantes $x$,$y$,$z$ sur \m{\P(V)} de telle sorte que \
\m{P_0=(1,0,0)}, \m{P_1=(0,1,0)}. Supposons que pour ces coordonn\'ees $\sigma$
soit d\'efinie par la matrice
\[
\left(\begin{array}{cccc}
0 & a & b & c+1 \\ -a & 0 & d-1 & e \\ -b & -d-1 & 0 & f \\ 1-c & -e & -f & 0
\end{array}\right).
\]
Alors en g\'en\'eral $D(\sigma)$ est constitu\'e de $L(\sigma)$ et des 
5 droites d\'efinies par l'\'equation
\[
(-ay-bx+(1-c)z)(ay^2+(c-d)yz+fz^2)(bx^2+(c+d)xz+ez^2) \ = \ 0.
\]

Supposons que $\sigma$ soit de type 2. Alors en g\'en\'eral $D(\sigma)$ est la 
r\'eunion de 4 droites de $\P(V)$. On choisit des
coordonn\'ees ind\'ependantes $x$,$y$,$z$ sur \m{\P(V)} de telle sorte que \
\m{P_0=(1,0,0)} et que $\ell$ soit la droite d'\'equation \ \m{z=0}. 
Supposons que pour ces coordonn\'ees $\sigma$
soit d\'efinie par la matrice
\[
\left(\begin{array}{cccc}
0 & a & b & c \\ -a & 0 & d & e-1 \\ -b & -d & 2 & f \\ -c & -e-1 & -f & 0
\end{array}\right).
\]
Alors en g\'en\'eral 
$D(\sigma)$ est constitu\'e de $L(\sigma)$ et des 3 droites d\'efinies par
l'\'equation
\[
(ax-dy-(e+1)z)(ay^2+byz+cz^2) \ = \ 0.
\]

Il existe n\'eanmoins des cas o\`u toutes les coniques de l'image de $\sigma$
sont d\'eg\'en\'er\'ees.
\end{sub}

\vskip 1.5cm

\begin{sub}\label{equiv}
{\bf \'Equivalence de congruences quadratiques planes}\rm

Soit $C$ une quadrique d'un hyperplan $\ell$ de \ \m{\P(V)=\P_2}, 
c'est-\`a-dire que $C$ est constitu\'ee de deux points distincts, ou est 
donn\'e par un seul point double \m{P_0} de $\ell$. Soit \m{W_0} l'espace des 
\'equations de coniques de \m{\P_2} contenant $C$.

\vskip 1cm

\begin{subsub}\label{theo4}{\bf Th\'eor\`eme : } Les fibres du morphisme
quotient d\'eduit de \m{\Phi_C}  
\[\Cong(C)/G_C\lra\Isom(W_0,W_0^*)/GL(W_0)\]
sont finies.\end{subsub}

\begin{proof} 
Pour \m{j=3,4} soit \m{I_j(W_0)} la sous-vari\'et\'e localement ferm\'ee
de \hfil\break\m{\Isom(W_0,W_0^*)} constitu\'ee des
$\phi$ tels que \m{Q(\phi)} soit de rang $j$. Ces sous-vari\'et\'es sont
\m{GL(W_0)}-invariantes.
On utilise les descriptions de \m{I_4(W_0))/GL(W_0)} et \m{I_3(W_0))/GL(W_0)}
d\'ecoulant du th\'eor\`eme
\ref{theo1}, dans les deux cas possibles pour $C$. 
Par exemple dans le premier cas, on consid\`ere des matrices
\[
\phi^{-1} \ = \ M \ = \ \left(\begin{array}{cccc}
0 & a & b & c+1 \\ -a & 0 & d-1 & e \\ -b & -d-1 & 0 & f \\
-c+1 & -e & -f & 0
\end{array}\right)
\]
(cf. la d\'emonstration de la proposition \ref{prop2}).
Tout \'el\'ement de \m{I_4(W_0)} est dans l'orbite d'une telle matrice,
relativement \`a la base indiqu\'ee dans la d\'emonstration de la proposition
\ref{prop2}. Le groupe \m{G_C} est constitu\'e des matrices \ 
\m{g = \left(\begin{array}{cccc}\alpha & 0 & u \\ 0 & \beta & v \\ 
0 & 0 & \gamma\end{array}\right)} \ avec \ \m{\alpha\beta\gamma\not= 0}. Alors
on a
\[
gM \ = \ \left(\begin{array}{cccc} 0 & a_1 & b_1 & c_1+\delta \\
-a_1 & 0 & d_1-\delta & e_1 \\ -b_1 & -d_1-\delta & 0 & f_1 \\
-c_1+\delta & -e_1 & -f_1 & 0
\end{array}\right)
\]
avec \ \m{\delta=\alpha\beta\gamma}. On peut supposer que \ \m{\delta=1} \ pour
que \m{gM} soit du m\^eme type que $M$, c'est-\`a-dire qu'on consid\`ere
l'action du sous-groupe \m{G_0} constitu\'e des $g$ tels que \
\m{\alpha\beta\gamma=1}. On a
\[
a_1=\alpha\gamma^2a, \ \ b_1=\beta\gamma^2b, \ \ c_1=c-\alpha\gamma av
-\beta\gamma bu, \ \ d_1=-\beta\gamma bu+\alpha\gamma av+d,
\]
\[
e_1=-\alpha\beta cu-\alpha\beta du+\alpha^2\beta e+\beta bu^2, \ \
f_1=\alpha\beta^2f+\alpha av^2-\alpha\beta cv+\alpha\beta dv.
\]
On est alors amen\'e \`a consid\'erer quatre ensembles \m{G_0}-invariants
de matrices $M$ :
\[
\km_1=\lbrace M; a\not=0, b\not=0\rbrace, \ \
\km_2=\lbrace M; a=0, b\not=0\rbrace,
\]
\[
\km_3=\lbrace M; a\not=0, b=0\rbrace, \ \
\km_4=\lbrace M; a=0, b=0\rbrace.
\]
Supposons que $M$ appartienne au premier ensemble. En faisant agir $G_0$ on se
ram\`ene au cas o\`u \ \m{a=b=1}, \m{e=f=0}. On obtient donc
\[
M \ = \ \left(\begin{array}{cccc}
0 & 1 & 1 & c+1 \\ -1 & 0 & d-1 & 0 \\ -1 & -d-1 & 0 & 0 \\
-c+1 & 0 & 0 & 0
\end{array}\right)
\]
On montre alors (en utilisant \'eventuellement un programme de calcul formel de
type Maple) que sauf si \ \m{c=d=0}, \m{\T(\phi)} est diagonalisable, et si \
\m{c=d=0}, on est dans le cas 2.5 de 2.2.2. Dans tous les cas
on peut mettre $c$ et $d$ de mani\`ere unique sous la forme
\[
c=\frac{1-\lambda}{1+\lambda}, \ \ d=\frac{\mu-1}{1+\mu}.
\]
(cela d\'ecoule du fait que $c$ et $d$ sont diff\'erents de $1$ et $-1$ car $M$
est non singuli\`ere). On montre alors que $\lambda$, \m{\dsp\frac{1}{\lambda}},
$\mu$ et  \m{\dsp\frac{1}{\mu}} sont des valeurs propres de \m{\T(\phi)}. Ceci
montre que l'application quotient \ \m{\km_1/G_0\to
I_4(W_0(P_0,P_1))/GL(W_0(P_0,P_1))} a des fibres finies. Les autres cas se
traitent de la m\^eme mani\`ere. 
\end{proof}



\end{sub}

\newpage

\section{Congruences quadratiques obtenues par translation d'une quadrique}

\begin{sub}\label{ExTr}\bf Exemples

\begin{subsub}Exemple 1\rm

Cet exemple appara\^{\i}t dans l'\'etude des fibr\'es exceptionnels sur \
\m{\P_1\times\P_1} (cf. \cite{ru1}, \cite{ru2}). Au point $P$ de coordonn\'ees
$x$, $y$, $z$ de \m{\P_2} on associe la conique \m{\sigma(x,y,z)} d'\'equation  
\[(zX-(x-z)Z)(zY-(y-z)Z)-z^2Z^2 \ = \ 0.\]
On obtient alors une congruence quadratique de type 1. Si on se limite au plan
r\'eel \ \m{\R^2\subset\P_2}, on associe au point \m{P=(x,y)} la translation de
vecteur \m{(x-1,y-1)} de l'hyperbole d'\'equation \ \m{XY=1} (voir la figure 1
ci-dessous). 

\null
\vspace{15cm}\hspace{-2cm}
\includegraphics{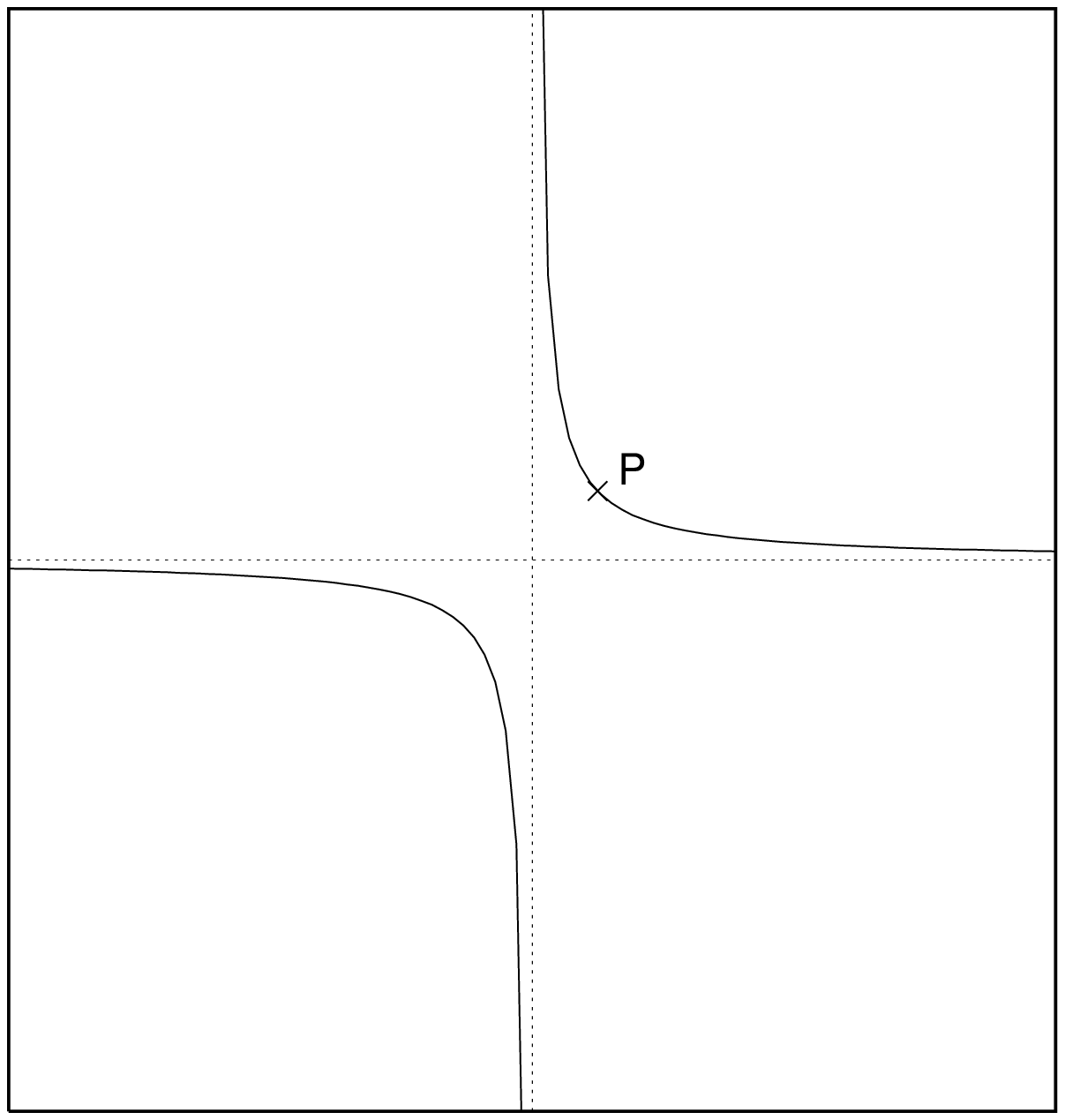}

\vspace{-2.5cm}

\begin{center}{\small Figure 1}
\end{center}

\end{subsub}

\newpage

\begin{subsub}\label{ExTr2}Exemple 2\rm

Cet exemple appara\^{\i}t dans l'\'etude des fibr\'es exceptionnels sur \
\m{\P_3} (cf. \cite{dr3}). Au point $P$ de coordonn\'ees
$x$, $y$, $z$ de \m{\P_2} on associe la conique \m{\sigma(x,y,z)} d'\'equation  
\[
-z^2XZ+(z^2+\frac{yz}{2})Y^2+(-\frac{y^2}{2}-2yz-\frac{4}{3}z^2)YZ+
(xz+y^2+\frac{4}{3}yz)Z^2 \ = \ 0.
\]
La congruence quadratique de type 2 ainsi obtenue devient plus explicite si on
fait le changement de coordonn\'ees (non lin\'eaire) suivant :
\[\theta : (X,Y,Z)\lra (YZ^2+\frac{1}{6}X^3,XZ^2,Z^3).\]
Alors \m{\theta^{-1}(\sigma(\theta(x,y,z)))} est la cubique d'\'equation
\[
-z^3YZ^2-\frac{1}{6}(zX-(x+2z)Z)^3+\frac{2}{3}z^2Z^2(zX-(x+2z)Z)+yz^2Z^3 \ = \
0.
\]
Si on se limite au plan
r\'eel \ \m{\R^2\subset\P_2}, on associe au point \m{P=(x,y)} la translation de
vecteur \m{(x+2,y)} de la cubique d'\'equation \ \m{Y \ = \ -\frac{1}{6}X^3-
\frac{2}{3}X} \ (voir la figure 2 ci-dessous). 

\null
\vspace{14cm}\hspace{-2cm}
\includegraphics{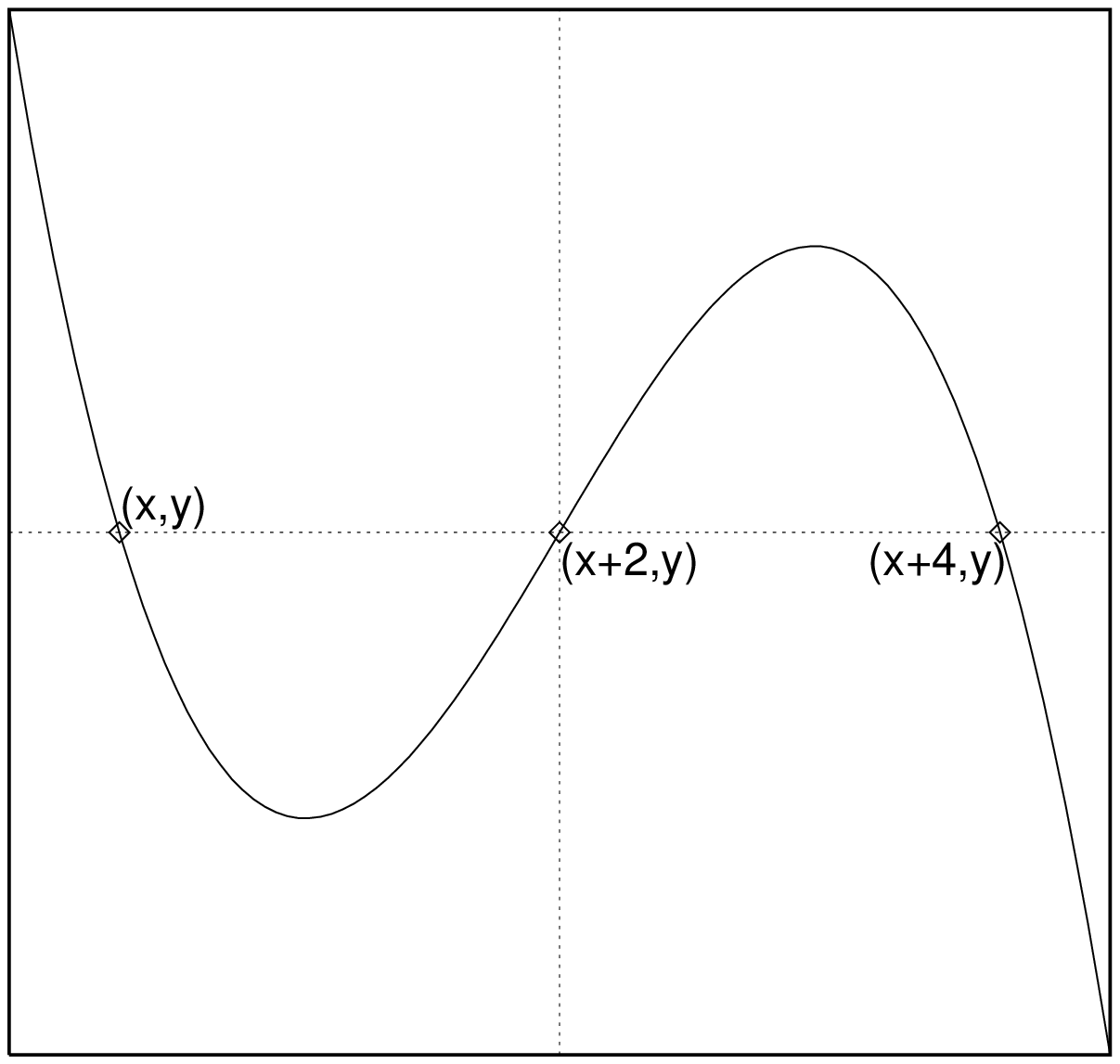}

\vspace{-2.5cm}

\begin{center}{\small Figure 2}
\end{center}

Dans ce syst\`eme de
coordonn\'ees, \m{\T(\sigma)} est la translation de vecteur \m{(4,0)} : dans la
figure pr\'ec\'edente, \m{(x+4,y)=\T(\sigma)(x,y)}.
\end{subsub}

\newpage
\begin{subsub}Exemple 3\rm

Soient $a$, $b$ des nombres complexes non nuls tels que \ \m{a+b\not=1}.
Au point $P$ de coordonn\'ees $x$, $y$, $z$ de \m{\P_2} on associe la conique 
\m{\sigma(x,y,z)} d'\'equation  
\[
z^2XY-axzYZ-byzXZ+(a+b-1)xyZ^2 \ = \ 0.
\]
Si $x$, $y$ et $z$ sont non nuls, c'est l'image inverse de la conique 
d'\'equation
\[ (X-aZ)(Y-bZ)-(1-a)(1-b)Z^2 \ = \ 0 \]
(qui ne d\'epend pas de \m{(x,y,z)}) par l'automorphisme de \m{\P_2}
\[ (X,Y,Z)\longmapsto (yzX,xzY,xyZ). \]
Si on se limite au plan r\'eel \ \m{\R^2\subset\P_2} et si $a$ et $b$ sont
r\'eels, on associe \`a \m{(x,y)} ($x$ et $y$ \'etant non nuls) l'image de 
l'hyperbole d'\'equation \ \m{(X-a)(Y-b)=(1-a)(1-b)} \ par l'automorphisme de 
\m{\R^2}
\[ (X,Y)\longmapsto (\frac{X}{x},\frac{Y}{y}). \]

\end{subsub}
\end{sub}

\vskip 1.5cm 

\begin{sub}\label{triv}\bf Sections de la quadrique universelle\rm

Soient $V$ un $\C$-espace vectoriel de dimension \m{n\geq 3}. 
Soit \ \m{\kq\subset\P(S^2V^*)\times\P(V)} \ la {\em quadrique universelle}. Il
est bien connu qu'il n'existe pas de section rationnelle \ \m{\P(S^2V^*)\to\kq}.
Si \ \m{W\subset S^2V^*} \ est un sous-espace vectoriel de dimension \m{n+1}, 
on note \m{\kq_W} la restriction de \m{\kq} \`a \m{\P(W)}. On va s'int\'eresser
\`a l'existence de sections rationnelles \ \m{\P(W)\to\kq_W}. 
On commencera par traiter le cas \ \m{n=3}, o\`u on sait caract\'eriser les
sous-espaces vectoriels $W$ tels qu'il existe une telle section rationnelle. On
traitera ensuite les cas o\`u \m{n>3} et o\`u $W$ est l'espace des \'equations 
de quadriques contenant une quadrique fixe d'un hyperplan de \m{\P(V)}. On
prouvera (par r\'ecurrence sur $n$) l'existence d'une section rationnelle de
\m{\kq_W}. 

On suppose maintenant que \ \m{n=3}. 

\vskip 1cm

\begin{subsub}\label{lemmepr3}{\bf Lemme : } Soient $U$,$S$,$T$ des
coordonn\'ees ind\'ependantes sur $V$. On consid\`ere le sous-espace vectoriel
$H$ de \m{S^2V^*} constitu\'e des coniques d'\'equation \ \m{aU^2+bS^2+cT^2=0}.
Alors il n'existe pas de section rationnelle \ \m{P(H)\lra\kq}.

\begin{proof}
Il faut montrer qu'il n'existe pas de polyn\^omes non nuls homog\`enes $P$,$Q$ 
et $R$ en $a$, $b$, $c$, de m\^emes degr\'es tels que \ \m{aP^2+bQ^2+cR^2=0}.
Supposons qu'il en existe. On peut supposer qu'ils sont premiers entre eux. En
se restreignant \`a l'hyperplan \m{a=0} on obtient l'\'equation \
\m{bQ^2=-cR^2} \ qui est impossible si $Q$ et $R$ sont non nuls sur \m{a=0}, car
dans le terme de gauche le facteur premier $c$ apparait avec un exposant pair et
dans le terme de droite il apparait avec un exposant impair. Il en d\'ecoule que
$Q$ et $R$ sont divisibles par $a$, et de m\^eme $b$ divise $P$ et $R$, et $c$
divise $P$ et $Q$. On peut donc \'ecrire
\[ P=bcP', \ \ Q=acQ', \ \ R=abR' \]
et on obtient l'\'equation
\[ bc{P'}^2+ac{Q'}^2+ab{R'}^2 \ = \ 0.  \]
En se restreignant \`a l'hyperplan \m{a=0} on voit encore que $P'$ est divisible
par $a$, et de m\^eme $Q'$ est divisible par $b$ et $R'$ par $c$. Finalement
$P$, $Q$ et $R$ sont divisibles par $abc$, contrairement \`a l'hypoth\`ese.
\end{proof}
\end{subsub}

\bigskip

Il existe 6 orbites distinctes de l'action de \m{PGL(V)} sur la grassmannienne
\m{Gr(4,S^2V^*)} des sous-espaces vectoriels de dimension 4 se \m{S^2V^*}. On le
voit plus facilement en consid\'erant la grassmannienne \m{Gr(2,S^2V)} des 
sous-espaces vectoriels de dimension 2 se \m{S^2V}. On trouve alors la liste
suivante, $V$ \'etant muni d'une base \m{(x,y,z)} :

\medskip

Cas 1 : \m{\pline{x^2,y^2}} -- espace des coniques de \m{\P(V^*)} qui sont des 
paires de droites contenant le point \m{(0,0,1)}, l'une \'etant l'image de
l'autre par l'involution de l'orthogonal de ce point \ \m{(x,y)\lra (x,-y)}. 
Le sous-espace de \m{S^2V^*} de dimension 4 correspondant est un espace de 
coniques de type 1.  

Cas 2 : \m{\pline{x^2,xy}} -- espace des coniques de \m{\P(V^*)} qui sont des
paires de droites contenant le point \m{(0,0,1)} et dont l'une est la droite
\m{x=0}. Le sous-espace de \m{S^2V^*} de dimension 4 correspondant est un 
espace de coniques de type 2.  

Cas 3 : \m{\pline{xy,xz}} -- espace des coniques de \m{\P(V^*)} qui sont des
paires de droites dont l'une est la droite \m{x=0} et l'autre passe par le point
\m{(0,0,1)}. Le sous-espace de \m{S^2V^*} de dimension 4 correspondant est un 
espace de coniques de type 3. 

Cas 4 : \m{\pline{xy,x(x+y+z)}} -- espace des coniques de \m{\P(V^*)} passant par
les points \m{(0,1,0)}, \m{(0,1,-1)}, \m{(1,0,0)} et \m{(1,0,-1)}. 
C'est le cas g\'en\'erique.

Cas 5 : \m{\pline{xy,z(x+y)}} -- espace des coniques de \m{\P(V^*)} passant par
\m{(1,0,0)}, \m{(0,1,0)} et \m{(0,0,1)} et dont la tangente en \m{(0,0,1)} est
la droite d'\'equation \ \m{x+y=0}.

Cas 6 : \m{\pline{xy,z^2}} -- espace des coniques de \m{\P(V^*)} contenant les
points \m{(1,0,0)}, \m{(0,1,0)}, dont la tangente en \m{(1,0,0)} est la droite
d'\'equation \m{y=0} et la tangente en \m{(0,1,0)} la droite d'\'equation
\m{x=0}.

\vskip 1cm

\begin{subsub}\label{prop3}{\bf Proposition : } Dans les cas 1,2 et 6 ci-dessus
il existe des sections rationnelles \ \m{\P(W)\to\Q_W} \ de la conique
universelle, et dans les autres cas il n'en existe pas.\rm

\begin{proof}
Les cas 1 et 2 sont imm\'ediats. Dans le cas 6, on consid\`ere la base 
\m{(XZ,YZ,X^2,Y^2)} de $W$ (\m{(X,Y,Z)} \'etant la base de \m{V^*} duale de
\m{(x,y,z)}). Le morphisme
\[\begin{array}{ccc}\P(W) & \lra & P(V) \\
aXZ+bYZ+cX^2+dY^2 & \longmapsto & (a+b,a+b,-c-d)
\end{array}
\]
d\'efinit une section rationnelle  \ \m{\P(W)\to\kq_W}. 

Il reste \`a montrer que dans les cas 3,4 et 5 il n'existe pas de section 
rationnelle. Une telle section induirait une section rationnelle de la conique
universelle sur tout 
hyperplan de \m{\P(W)}. Notons que dans les cas 3,4 et 5 il existe toujours des
coordonn\'ees $U$, $S$, $T$ sur $V$ telles que $P(W)$ contienne toutes les
coniques d'\'equation \ \m{aU^2+bS^2+cT^2=0} ($a$, $b$ et $c$ parcourant $\C$).  
Le r\'esultat est donc une cons\'equence du lemme \ref{lemmepr3}.
\end{proof}
\end{subsub}

\bigskip

Soit \ \m{W\subset S^2V} \ un espace de coniques de type 1 ou 2,
de dimension 4. Rappelons (cf. chapitre 4) qu'on dit que $W$ est un espace de
coniques de type 1 (resp. 2) s'il existe deux points distincts \m{P_0}, \m{P_1}
de \m{\P(V)} (resp. un point \m{P_0} et une droite \m{\ell} de \m{\P(V)} passant
par \m{P_0}) tels que $W$ soit l'espace des \'equations de coniques passant par
\m{P_0} et \m{P_1} (resp. passant par \m{P_0} et tangentes \`a \m{\ell} en
\m{P_0}).
On note \m{W^0} l'ouvert de $W$ constitu\'e des coniques lisses.

\vskip 1cm

\begin{subsub}\label{prop4}{\bf Proposition : }
La restriction \m{\kq_{W^0}} de \m{\kq} \`a \m{\P(W^0)} est un fibr\'e en
espaces projectifs trivial.\rm 

\begin{proof}
Pour le type 1 on consid\`ere une conique d'\'equation 
\[uXY+vYZ+wXZ+tZ^2 \ = \ 0.\]
La condition de lissit\'e pour cette conique est \ \m{u(vw-ut)\not=0}. 
La matrice 
\[
\left(\begin{array}{ccc}
u(vw-ut) & 0 & v(vw-ut) \\ 0 & u & w \\ 0 & 0 & vw-ut
\end{array}\right)
\]
d\'efinit un automorphisme de \m{\P(V)} qui induit un isomorphisme entre la
conique pr\'ec\'edente et la conique d'\'equation \ \m{XY-Z^2=0}. La proposition
\ref{prop4} en d\'ecoule imm\'ediatement. Le cas du type 2 est analogue.
\end{proof}
\end{subsub}

\vskip 1cm

\begin{subsub}\label{casgen}{\bf Le cas \ \m{n>3 \ }. }\rm 
On suppose ici que $V$ est de dimension \ \m{n>3}. Le r\'esultat suivant 
g\'en\'eralise une partie des propositions \ref{prop3} et \ref{prop4} :

\vskip 1cm

\begin{subsub}\label{prop8}
{\bf : Proposition : } Soient $C$ une quadrique d'un hyperplan de
\m{\P_{n-1}}, $W$ le sous-espace vectoriel de dimension \m{n+1} de \m{S^2V^*}
constitu\'e des \'equations des quadriques contenant $C$, et \m{W^0} l'ouvert de
$W$ correspondant aux quadriques de rang maximal. Alors il existe des sections
de \m{\kq_{W^0}}. Ce dernier est un fibr\'e en quadriques trivial sur
\m{\P(W^0)}.
\end{subsub}

\begin{proof}
Il suffit de montrer qu'il existe un morphisme
\[\Phi : \P(W^0)\lra\Aut(\P_{n-1})\]
et une quadrique \m{Q_0} de \m{\P_{n-1}} tels que pour tout \m{Q\in\P(W^0)},
$Q$ (vu comme quadrique) soit l'image r\'eciproque de \m{Q_0} par \m{\Phi(Q)}.
On se ram\`eme ais\'ement au cas o\`u les quadriques de \m{\P(W_0)} sont de rang
$n$.

Si \m{n=3}, en modifiant la proposition \ref{prop4} on voit que la quadrique
d'\'equation 
\[ X_0^2+X_1^2+X_2(aX_0+bX_1+cX_2) \ = \ 0
\]
est de rang 3 si et seulement si \ \m{\delta=a^2+b^2-4c} \ est non nul, et que
dans ce cas cette quadrique est l'image r\'eciproque de la quadrique 
d'\'equation \ \m{X_0^2+X_1^2-X_2^2 = 0} \ par l'automorphisme de \m{\P_2} 
donn\'e par la matrice
\[
\left(\begin {array}{ccc} \frac{1}{2}(\delta+1) & \frac{1}{2}\,
i\left (\delta-1\right ) &
\frac{1}{4}(\delta(a-ib)+a+ib) 
\\\noalign{\medskip}\frac{1}{2}\,i\left (\delta-1\right )&
-\frac{1}{2}(\delta+1)&\frac{1}{4}\,i\left (\delta(a-ib)-a-
ib\right )\\\noalign{\medskip}0&0&\frac{1}{2}\,\delta\end {array}
\right )
. \]
Cette matrice permet ais\'ement de construire $\Phi$ dans ce cas.

Supposons le r\'esultat prouv\'e pout \m{n-1}. La quadrique d'\'equation
\[
X_0^2+\cdots+X_{n-2}^2+X_{n-1}(\alpha_0X_0+\cdots+\alpha_{n-1}X_{n-1}) \ = \ 0
\]
est de rang $n$ si et seulement si
\m{\alpha_0^2+\cdots+\alpha_{n-2}^2-4\alpha_{n-1}\not= 0.}
Supposons que ce soit le cas.
Soit \m{\gamma} un nombre complexe tel que \
\m{\alpha_0^2+\cdots+\alpha_{n-3}^2-4\gamma\not= 0}. Ceci implique que la
quadrique d'\'equation \ 
\m{ X_0^2+\cdots+X_{n-3}^2+X_{n-1}(\alpha_0X_0+\cdots+\alpha_{n-3}X_{n-3}
+\gamma X_{n-1})=0} \ est de rang \m{n-1}. D'apr\`es l'hypoth\`ese de
r\'ecurrence, il existe une matrice \
\m{M=(a_{ij})_{0\leq i,j \leq n-3}} \ non singuli\`ere, et des scalaires
\m{\mu_0,\ldots,\mu_{n-3},\beta}, d\'ependant
alg\'ebriquement de \m{\alpha_0,\ldots,\alpha_{n-3},\gamma} et tels que
\[ 
\sigg_{i=0}^{n-3}(\sigg_{j=0}^{n-3}a_{ij}X_j+\mu_iX_{n-1})^2 \ - 
\beta^2X_{n-1}^2 \ = \ \ \ \ \ \ \ \ \ \ \ \ \ \ \ \ \ \ \ \ \ \ \ \ \ \]
\[ \ \ \ \ \ \ \ \ \ \ \ \ \ \ \ \ \ \ \ \ \ \ \ \ \ \ \ \ \ \
X_0^2+\cdots+X_{n-3}^2+X_{n-1}(\alpha_0X_0+\cdots+\alpha_{n-3}X_{n-3}+\gamma
X_{n-1}).
\]
Soient \m{\mu_{n-2}, \mu_{n-1}} des nombres complexes. Alors on a
\[
\sigg_{i=0}^{n-3}(\sigg_{j=0}^{n-3}a_{ij}X_j+\mu_iX_{n-1})^2 \ +
(X_{n-2}+\mu_{n-2}X_{n-1})^2-\mu_{n-1}^2X_{n-1}^2 \ = \
\ \ \ \ \ \ \ \ \ \ \ \ \ \ \ \ \ \ \ \ \ \ \ \ \ \ \]
\[
X_0^2+\cdots+X_{n-2}^2+(\beta^2-\mu_{n-1}^2+\mu_{n-2}^2+\gamma)X_{n-1}^2 +
\ \ \ \ \ \ \ \ \ \ \ \ \ \ \ \ \ \ \ \ \ \ \ \ \ \ \] 
\[\ \ \ \ \ \ \ \ \ \ \ \ \ \ \ \ \ \ \ \ \ \ \ \ \ \ 
X_{n-1}(\alpha_0X_0+\cdots+\alpha_{n-3}X_{n-3}+2\mu_{n-2}X_{n-2}).
\]
On prend maintenant 
\[ \gamma \ = \alpha_{n-1} - \frac{\alpha_{n-2}^2}{4}. \] 
Remarquons qu'on a bien
\ \m{\alpha_0^2+\cdots+\alpha_{n-3}^2-4\gamma\not= 0}. On prend ensuite \
\m{\mu_{n-1}=\beta}. On a alors
\[
\sigg_{i=0}^{n-3}(\sigg_{j=0}^{n-3}a_{ij}X_j+\mu_iX_{n-1})^2 \ +
(X_{n-2}+\mu_{n-2}X_{n-1})^2-\mu_{n-1}^2X_{n-1}^2 \ = \
\ \ \ \ \ \ \ \ \ \ \ \ \ \ \ \ \ \ \ \ \ \ \ \ \ \ \]
\[ \ \ \ \ \ \ \ \ \ \ \ \ \ \ \ \ \ \ \ \ \ \ \ \ \ \ 
X_0^2+\cdots+X_{n-2}^2+X_{n-1}(\alpha_0X_0+\cdots+\alpha_{n-1}X_{n-1}),
\]
ce qui d\'efinit le morphisme voulu $\Phi$.
\end{proof}
\end{subsub}
\end{sub}

\vskip 1.5cm

\begin{sub}\label{xx0} \bf
Congruences quadratiques obtenues par translations d'une
quadrique\rm

Soit $V$ un \m{\C}-espace vectoriel de dimension $n$.

\bigskip

\begin{subsub}\label{th0b}{\bf Th\'eor\`eme : }
Soient $C$ une quadrique d'un hyperplan de
\m{ \P(V)}, $W$ le sous-espace vectoriel de dimension \m{n+1} de \m{S^2V^*}
constitu\'e des \'equations des quadriques contenant $C$ et
\[\sigma :\P(V)\lra\P(W)\]
une congruence quadratique dont l'image contient au moins une quadrique de rang
maximal de \m{\P(W)}. Alors il existe un morphisme rationnel
\[ R : \P(V)\lra PGL(V) \]
et une quadrique \m{C_0} de \m{\P(V)} tels que pour un point g\'en\'eral $P$ de
\m{\P_{n-1}}, on ait
\[ R(P)^{-1}(C_0) \ = \ \sigma(P). \]
\end{subsub}

\begin{proof}
D\'ecoule imm\'ediatement de la proposition \ref{prop8}.
\end{proof}

\vskip 1cm

\begin{subsub}\label{remx}{\bf Remarque : }\rm
L'image de $\sigma$ contient toujours une quadrique de \m{\P(W)} de rang maximal
si le rang de $C$ n'est pas maximal. Dans le cas contraire il se peut que toute
quadrique de l'image de $\sigma$ ne soit pas de rang maximal.
\end{subsub}

\vskip 1cm

\begin{subsub}\label{cas_plan}{\bf Cas des congruences quadratiques planes. }\rm
On va donner une formule explicite pour $R$ dans le cas \ \m{n=3}.
Soit $\sigma$ une congruence quadratique de type 1. On choisit des
coordonn\'ees ind\'ependantes $X$,$Y$,$Z$ sur \m{\P(V)} de telle sorte que \
\m{P_0=(1,0,0)}, \m{P_1=(0,1,0)}. Supposons que pour ces coordonn\'ees $\sigma$
soit d\'efinie par la matrice
\[
\left(\begin{array}{cccc}
0 & a & b & c+1 \\ -a & 0 & d-1 & e \\ -b & -d-1 & 0 & f \\ 1-c & -e & -f & 0
\end{array}\right).
\]
Posons
\[
\begin{array}{cc}
\phi_{XY}(x,y,z)=-ayz-bxz+(1-c)z^2, & \phi_{YZ}(x,y,z)=axy-(d+1)xz-ez^2, \\
\phi_{XZ}(x,y,z)=bxy+(d-1)yz-fz^2, & \phi_{Z^2}(x,y,z)=(c+1)xy+eyz+fxz.
\end{array}
\]
On a alors
\[
\sigma(x,y,z) \ = \ \phi_{XY}(x,y,z)XY+\phi_{YZ}(x,y,z)YZ+\phi_{XZ}(x,y,z)XZ+
\phi_{Z^2}(x,y,z)Z^2.
\]
Soit
\[
\Delta \ = \ \phi_{YZ}\phi_{XZ}-\phi_{XY}\phi_{Z^2}.
\]
Soit R le morphisme rationnel \ \m{\P(V)\lra\P(\End(V))} d\'efini par la matrice
\[
\left(\begin{array}{ccc}
\phi_{XY}\Delta & 0 & \phi_{YZ}\Delta \\ 0 & z^4\phi_{XY} & z^4\phi_{XZ} \\
0 & 0 & z^2\Delta
\end{array}\right).
\]
Alors on v\'erifie imm\'ediatement que $R$ poss\`ede la propri\'et\'e du 
th\'eor\`eme \ref{th0b}.

Le cas des congruences quadratiques de type 2 est beaucoup plus simple.
On choisit des
coordonn\'ees ind\'ependantes $X$,$Y$,$Z$ sur \m{\P(V)} de telle sorte que \
\m{P_0=(1,0,0)} et que l'\'equation de $\ell$ soit \ \m{Z=0}. 
Supposons que pour ces coordonn\'ees $\sigma$ soit d\'efinie par la matrice
\[
\left(\begin{array}{cccc}
0 & a & b & c \\ -a & 0 & d & e-1 \\ -b & -d & 2 & f \\ -c & -e-1 & -f & 0
\end{array}\right).
\]
Posons
\[
\begin{array}{cc}
\phi_{XZ}(x,y,z)=-ay^2-byz-cz^2, & \phi_{Y^2}(x,y,z)=axz-dyz-(e+1)z^2, \\
\phi_{YZ}(x,y,z)=bxz+dy^2+2yz-fz^2, & \phi_{Z^2}(x,y,z)=cxz+(e-1)y^2+fyz.
\end{array}
\]
On a alors
\[
\sigma(x,y,z) \ = \ \phi_{XZ}(x,y,z)XZ+\phi_{Y^2}(x,y,z)Y^2+\phi_{YZ}(x,y,z)YZ+
\phi_{Z^2}(x,y,z)Z^2.
\]
On d\'efinit alors $R$ par la matrice
\[
\left(\begin{array}{ccc}
\phi_{XZ} & \phi_{YZ} & \phi_{Z^2} \\ 0 & 0 & -\phi_{Y^2} \\ 
0 & \phi_{Y^2} & 0
\end{array}\right).
\]

\vskip 1cm

{\bf Remarques : }\rm
\begin{enumerate}
\item Dans ce qui pr\'ec\`ede les polyn\^omes 
d\'efinissant $R$ pour les congruences quadratiques de type 1 
sont de degr\'e 6 (cf. \ref{propx}).
\item Le lieu des points de $\P(V)$ o\`u les morphismes $R$ pr\'ec\'edents
ne sont pas d\'efinis est exactement 
$D(\sigma)$ (cf. \para \ref{degen}).
\end{enumerate}
\end{subsub}

\vskip 1cm

\begin{subsub}\label{propx}{\bf Proposition : } Soit
\[\sigma : \P_2=\P(V)\lra\P_3=\P(W)\]
une congruence
quadratique plane telle que \ \m{L(\sigma)=H(\sigma)}. Alors il existe un
morphisme rationnel 
\[R:\P(V)\lra PGL(V)\]
d\'efini par des formes quadratiques et une conique \m{C_0} de \m{\P(V)}, tels 
que pour un point g\'en\'eral $P$ de \m{\P(V)} on ait \ \m{\sigma(P)=R^{-1}(Q)}.
\end{subsub}

\begin{proof}
Soient $\alpha_0$, $\alpha_1$, $\beta_0$, $\beta_1$, $e$,$f$ des nombres 
complexes tels que 
\[\alpha_0\not=\alpha_1, \ \
\alpha_0^2-4\alpha_0\alpha_1+\alpha_1^2-\beta_0-\beta_1=2.\]
Soient
\[
q_0 \ = \ Z((\beta_0+\alpha_1^2)Y-fZ), \ \ \ q_1 \ = \ Z((\beta_1+\alpha_0^2)X
-eZ),
\]
\[
s_0 \ = \ (\beta_0\beta_1-\alpha_0^2\alpha_1^2)XY-\beta_0eYZ-\beta_1fXZ+efZ^2,
\]
\[
s_1 = (-\beta_0-\beta_1-\alpha_0^2-\alpha_1^2)XY+eYZ+fXZ, 
\]
\[
s_2 \ = \
-(\alpha_0\alpha_1^2+\alpha_0^2\alpha_1+\alpha_0\beta_0+\alpha_1\beta_1)XY
+\alpha_1eYZ+\alpha_0fXZ.
\]
On d\'efinit $R$ par 
\[
R(X,Y,Z) \ = \ \left(\begin{array}{ccc}
\alpha_0^2q_0 & \alpha_1^2q_1 & s_0 \\ q_0 & q_1 & s_1 \\
\alpha_0q_0 & \alpha_1q_1 & s_2
\end{array}\right) .
\]
se \m{C_0} est la conique d\'equation \ \m{XY-Z^2=0}.
Dans ce cas $\sigma$ est d\'efini par la matrice 
\[
\left(\begin{array}{cccc}
0 & 0 & 0 & c+1 \\ 0 & 0 & d-1 & e \\ 0 & -d-1 & 0 & f \\ -c+1 & -e & -f & 0
\end{array}\right)
\]
avec
\[
c \ = \ 1-(\alpha_0-\alpha_1)^2, \ \
d \ = \ 1+\beta_0-\alpha_0^2+2\alpha_0\alpha_1.
\]
Pour le voir on calcule le produit
\[R(X,Y,Z)\left(\begin{array}{c}x \\ y \\ z \end{array}\right) \ = \
\left(\begin{array}{c} X_0 \\ Y_0 \\ Z_0 \end{array}\right). \]
L'\'equation de \m{\sigma(X,Y,Z)} n'est pas exactement \ \m{X_0Y_0-Z_0^2=0}.
On constate en fait que  \ \m{\psi(x,y,z)=X_0Y_0-Z_0^2} \ est divisible par
\m{(\beta_0+\alpha_1^2)Y-fZ} et \m{(\beta_1+\alpha_0^2)X-eZ} :
\[\psi(x,y,z) \ = \ ((\beta_0+\alpha_1^2)Y-fZ)((\beta_1+\alpha_0^2)X-eZ)
q(x,y,z),\]
et l'\'equation de \m{\sigma(X,Y,Z)} est \ \m{q(x,y,z)=0}.
\end{proof}

\vskip 1cm

\begin{subsub}{\bf Remarque : }\rm Dans certains cas il est possible que $R$
puisse \^etre d\'efini par des formes lin\'eaires.
Soient $\nu$, $\alpha_0$, $\alpha_1$,$w$ des nombres complexes tels que 
\[
\alpha_0\not=\alpha_1, \ \ \nu=\frac{1}{(\alpha_0-\alpha_1)^2}, \ \
w\not=0.
\]
On d\'efinit $R$ par
\[
R(X,Y,Z) \ = \ \left(\begin{array}{ccc}
\nu\alpha_0^2Z & \alpha_1^2Z & -\nu\alpha_0^2X-\alpha_1^2Y+wZ \\ 
\nu Z & Z & -\nu X -Y \\
\nu\alpha_0Z & \alpha_1Z & -\nu\alpha_0X-\alpha_1Y
\end{array}\right)
\]
(\m{C_0} \'etant la conique d\'equation \ \m{XY-Z^2=0}).
Dans ce cas $\sigma$ est d\'efini par la matrice
\[
\left(\begin{array}{cccc}
0 & 0 & 0 & 1 \\ 0 & 0 & -1 & -w \\
0 & -1 & 0 & -\frac{w}{(\alpha_0-\alpha_1)^2} \\
1 & w & \frac{w}{(\alpha_0-\alpha_1)^2} & 0
\end{array}\right)\]
\end{subsub}

\vskip 1cm

\begin{subsub}{\bf Conjecture : } Dans le cas g\'en\'eral, une congruence
quadratique normale \hfil\break
\m{\sigma : \P_{n-1}\lra\P(W)} \ provient d'un morphisme
\ \m{R:\P_{n-1}\lra PGL(n)} \ d\'efini par des formes quadratiques si et
seulement si on a \ \m{L(\sigma)=H(\sigma)}.
\end{subsub}
\end{sub}

\newpage

\begin{sub}\bf Congruences quadratiques planes
obtenues par fixation d'une tangente\rm

\begin{subsub}Exemple\rm

Au point $P$ de coordonn\'ees
$x$, $y$, $z$ de \m{\P_2} on associe la conique \m{\sigma(x,y,z)} d'\'equation  
\[
yzXZ-\frac{z^2}{2}Y^2+(yz-xz)YZ-\frac{y^2}{2}Z^2 \ = \ 0.
\]
On obtient ainsi une congruence quadratique de type 2. On a
\[ \T(\sigma)(x,y,z) \ = \ (x-2z,y,-z). \]
Cette congruence quadratique peut s'interpr\'eter de la fa\c con suivante : 
\m{\sigma(x,y,z)} est l'unique conique passant par les points \m{(1,0,0)},
\m{(x,y,z)}, \m{\T(\sigma)(x,y,z)}, dont la tangente en \m{(1,0,0)} est la
droite  \m{Z=0}, et la tangente en \m{(x,y,z)} la droite d'\'equation \
\m{Xy-(x+y)Y=0}. Plus concr\`etement, on se limite \`a \ \m{\R^2\subset\P_2},
les coordonn\'ees r\'eelles \'etant $X$ et $Z$. Au point $P$ de \m{\R^2} de
coordonn\`ees $x$, $z$ on associe l'unique hyperbole dont l'une des asymptotes
est la droite d'\'equation \ \m{Z=0}, qui contient $P$ et \m{(x-2,-z)}, et dont
la tangente en $P$ est la droite d'\'equation \ \m{X=x} (cf. figure 3
ci-dessous).

\null
\vspace{14cm}\hspace{-2cm}
\includegraphics{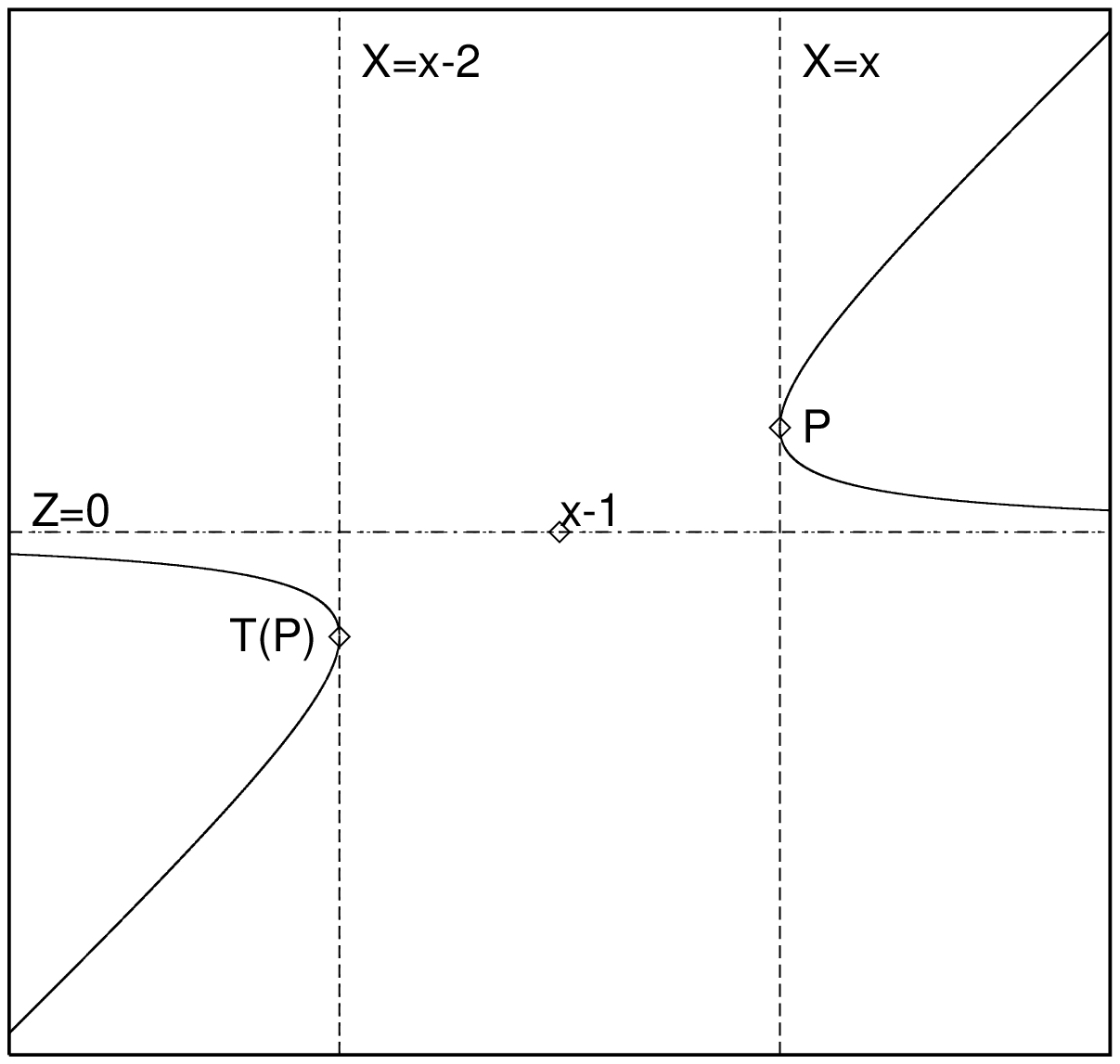}

\vspace{-2.5cm}

\begin{center}{\small Figure 3}
\end{center}
\end{subsub}

\newpage

\begin{subsub}Cas g\'en\'eral\rm

On ne s'\'etendra pas sur le sujet. On peut compl\`etement d\'ecrire une
congruence quadratique $\sigma$ en sp\'ecifiant pour tout point $P$ de \m{\P_2},
un point de \m{\sigma(P)} autre que $P$ (par exemple \m{\T(\sigma)(P)}) 
et la tangente \`a \m{\sigma(P)} en $P$. On peut obtenir de telles descriptions
en utilisant le th\'eor\`eme \ref{th0b}.
\end{subsub}
\end{sub}

\vskip 2.5cm

\section{Liste des types de congruences quadratiques}

On donne ci-dessous la liste des valeurs possibles des orbites de $\sigma$ 
correspondant aux translations \'enonc\'ees en \ref{theo1}, pour les
deux premi\`eres situations de la proposition \ref{prop1}. Pour les
r\'ealisations g\'eom\'etriques par des coniques de type 1, on consid\`ere les
orbites par le groupe laissant \m{\lbrace P_0,P_1\rbrace} invariant.

\vskip 1.5cm

\begin{sub}\label{liste1}\bf 
R\'ealisation g\'eom\'etrique par des coniques de type
1\rm

{\bf Cas 1.1.} On obtient 3 orbites principales. On obtient d'autres orbites en
rempla\c cant $\lambda$ par \m{1/\lambda} ou $\mu$ par \m{1/\mu}. 
Les matrices \ \m{\sigma^*=\phi^{-1}} \ sont les suivantes :

\[
\phi^{-1}_{1.1a} \ = \ \left(\begin{array}{cccc}
0 & 1 & 1 & \frac{2}{\lambda+1} \\ -1 & 0 & -\frac{2}{\mu+1} & 0 \\
-1 & -\frac{2\mu}{\mu+1} & 0 & 0 \\ \frac{2\lambda}{\lambda+1} & 0 & 0 & 0
\end{array}\right), \ \ \ \
\phi^{-1}_{1.1b} \ = \ \left(\begin{array}{cccc}
0 & 0 & 0 & \frac{2}{\lambda+1} \\ 0 & 0 & -\frac{2\mu}{\mu+1} & 0 \\
0 & -\frac{2}{\mu+1} & 0 & 0 \\ \frac{2\lambda}{\lambda+1} & 0 & 0 & 0
\end{array}\right),
\]
\[
\phi^{-1}_{1.1c} \ = \ \left(\begin{array}{cccc}
0 & 1 & 0 & 2 \\ -1 & 0 & -\frac{2(\lambda+\mu+2)}{(\lambda+1)(\mu+1)} & 0 \\
0 & -\frac{2(\lambda\mu-1)}{(\lambda+1)(\mu+1)} & 0 & 
\frac{4}{(\lambda+1)(\mu+1)} \\ 
0 & 0 & -\frac{4}{(\lambda+1)(\mu+1)} & 0 
\end{array}\right).
\]
\end{sub}

\vskip 1cm

{\bf Cas 1.2.} On obtient 3 orbites principales. On obtient d'autres orbites en
rempla\c cant $\lambda$ par \m{1/\lambda}. Les matrices \
\m{\sigma^*=\phi^{-1}} \ sont les suivantes :

\[
\phi^{-1}_{1.2a} \ = \ \left(\begin{array}{cccc}
0 & 1 & 0 & 2 \\ -1 & 0 & -2 & 0 \\ 
0 & 0 & 0 & \frac{4\lambda}{(1+\lambda)^2} \\ 
0 & 0 & \frac{-4\lambda}{(1+\lambda)^2} & 0
\end{array}\right), \ \ \ \
\phi^{-1} \ = \ \left(\begin{array}{cccc}
0 & 0 & 0 & \frac{2}{\lambda+1} \\ 0 & 0 & -\frac{2\lambda}{\lambda+1} & 0 \\
0 & -\frac{2}{\lambda+1} & 0 & 0 \\ \frac{2\lambda}{\lambda+1} & 0 & 0 & 0
\end{array}\right), 
\]
\[
\phi^{-1}_{1.2c} \ = \ \left(\begin{array}{cccc}
0 & 0 & 0 & \frac{2}{\lambda+1} \\ 
0 & 0 & -\frac{2\lambda}{\lambda+1} & 1 \\
0 & -\frac{2}{\lambda+1} & 0 & 0 \\ 
\frac{2\lambda}{\lambda+1} & -1 & 0 & 0 
\end{array}\right).
\]

\vskip 1cm

{\bf Cas 1.5.} On obtient 9 orbites principales. On obtient d'autres orbites en
rempla\c cant $\lambda$ par \m{1/\lambda}. Les matrices \
\m{\sigma^*=\phi^{-1}} \ sont les suivantes :
\[
\phi^{-1}_{1.5a} \ = \ \left(\begin{array}{cccc}
0 & 1 & 1 & 1 \\ -1 & 0 & -\frac{2}{\lambda+1} & 0 \\ 
-1 & -\frac{2\lambda}{\lambda+1} & 0 & 0 \\ 
1 & 0 & 0 & 0
\end{array}\right), \ \ \ \
\phi^{-1}_{1.5b} \ = \ \left(\begin{array}{cccc}
0 & 1 & 1 & \frac{2}{\lambda+1} \\ -1 & 0 & -1 & 0 \\
-1 & -1 & 0 & 0 \\ \frac{2\lambda}{\lambda+1} & 0 & 0 & 0
\end{array}\right), 
\]
\[
\phi^{-1}_{1.5c} \ = \ \left(\begin{array}{cccc}
0 & 1 & 0 & 2 \\ -1 & 0 & -\frac{\lambda+3}{\lambda+1} & 0 \\
0 & -\frac{\lambda-1}{\lambda+1} & 0 & \frac{2}{\lambda+1} \\ 
0 & 0 & -\frac{2}{\lambda+1} & 0 
\end{array}\right), \ \ \ \
\phi^{-1}_{1.5d} \ = \ \left(\begin{array}{cccc}
0 & 0 & 0 & \frac{2}{\lambda+1} \\ 0 & 0 & -1 & 0 \\
0 & -1 & 0 & 0 \\ \frac{2\lambda}{\lambda+1} & 0 & 0 & 0 
\end{array}\right),
\]
\[
\phi^{-1}_{1.5e} \ = \ \left(\begin{array}{cccc}
0 & 0 & 0 & \frac{2}{\lambda+1} \\ 0 & 0 & -1 & 1 \\
0 & -1 & 0 & 0 \\ \frac{2\lambda}{\lambda+1} & -1 & 0 & 0 
\end{array}\right), \ \ \ \
\phi^{-1}_{1.5f} \ = \ \left(\begin{array}{cccc}
0 & 0 & 0 & \frac{2}{\lambda+1} \\ 0 & 0 & -1 & 1 \\
0 & -1 & 0 & 1 \\ \frac{2\lambda}{\lambda+1} & -1 & -1 & 0 
\end{array}\right),
\]
\[
\phi^{-1}_{1.5g} \ = \ \left(\begin{array}{cccc}
0 & 0 & 0 & 1 \\ 0 & 0 & -\frac{2\lambda}{\lambda+1} & 0 \\
0 & -\frac{2}{\lambda+1} & 0 & 0 \\ 1 & 0 & 0 & 0 
\end{array}\right), \ \ \ \
\phi^{-1}_{1.5h} \ = \ \left(\begin{array}{cccc}
0 & 0 & 0 & 1 \\ 0 & 0 & -\frac{2\lambda}{\lambda+1} & 1 \\
0 & -\frac{2}{\lambda+1} & 0 & 0 \\ 1 & -1 & 0 & 0 
\end{array}\right),
\]
\[
\phi^{-1}_{1.5i} \ = \ \left(\begin{array}{cccc}
0 & 0 & 0 & 1 \\ 0 & 0 & -\frac{2\lambda}{\lambda+1} & 1 \\
0 & -\frac{2}{\lambda+1} & 0 & 1 \\ 1 & -1 & -1 & 0 
\end{array}\right).
\]

\vskip 1cm

{\bf Cas 2.3.} On obtient 5 orbites principales. On obtient d'autres orbites en
rempla\c cant $\lambda$ par \m{1/\lambda}. Les matrices \
\m{\sigma^*=\phi^{-1}} \ sont les suivantes :
\[
\phi^{-1}_{2.3a} \ = \ \left(\begin{array}{cccc}
0 & 1 & 1 & \frac{2}{\lambda+1} \\ -1 & 0 & -\frac{2}{\lambda+1} & 0 \\ 
-1 & -\frac{2\lambda}{\lambda+1} & 0 & 0 \\ 
\frac{2\lambda}{\lambda+1} & 0 & 0 & 0
\end{array}\right), \ \ \ \
\phi^{-1}_{2.3b} \ = \ \left(\begin{array}{cccc}
0 & 1 & 0 & 2 \\ -1 & 0 & -\frac{4}{\lambda+1} & 0 \\
0 & -\frac{2(\lambda-1)}{\lambda+1} & 0 & \frac{4}{(1+\lambda)^2} \\ 
0 & 0 & -\frac{4}{(1+\lambda)^2} & 0
\end{array}\right), 
\]
\[
\phi^{-1}_{2.3c} \ = \ \left(\begin{array}{cccc}
0 & 0 & 1 & 2 \\ 0 & 0 & 0 & \frac{4\lambda}{(1+\lambda)^2} \\
-1 & -2 & 0 & 1 \\ 0 & -\frac{4\lambda}{(1+\lambda)^2} & -1 & 0 
\end{array}\right), \ \ \ \
\phi^{-1}_{2.3d} \ = \ \left(\begin{array}{cccc}
0 & 0 & 0 & \frac{2}{\lambda+1} \\ 0 & 0 & -\frac{2\lambda}{\lambda+1} & 0 \\
0 & -\frac{2}{\lambda+1} & 0 & 1 \\ \frac{2\lambda}{\lambda+1} & 0 & -1 & 0 
\end{array}\right),
\]
\[
\phi^{-1}_{2.3e} \ = \ \left(\begin{array}{cccc}
0 & 0 & 0 & \frac{2}{\lambda+1} \\ 0 & 0 & -\frac{2\lambda}{\lambda+1} & 1 \\
0 & -\frac{2}{\lambda+1} & 0 & 1 \\ \frac{2\lambda}{\lambda+1} & -1 & -1 & 0 
\end{array}\right).
\]

\vskip 1cm

{\bf Cas 2.4.} On obtient 2 orbites. Les matrices \
\m{\sigma^*=\phi^{-1}} \ sont les suivantes :
\[
\phi^{-1}_{2.4a} \ = \ \left(\begin{array}{cccc}
0 & 1 & 0 & 2 \\ -1 & 0 & -2 & 0 \\ 0 & 0 & 0 & 1 \\ 0 & 0 & -1 & 0
\end{array}\right), \ \ \ \
\phi^{-1}_{2.4b} \ = \ \left(\begin{array}{cccc}
0 & 0 & 0 & 1 \\ 0 & 0 & -1 & 1 \\ 0 & -1 & 0 & 0 \\ 1 & -1 & 0 & 0
\end{array}\right).
\]

\vskip 1cm

{\bf Cas 2.5.} On obtient 3 orbites. Les matrices \
\m{\sigma^*=\phi^{-1}} \ sont les suivantes :
\[
\phi^{-1}_{2.5a} \ = \ \left(\begin{array}{cccc}
0 & 1 & 1 & 1 \\ -1 & 0 & -1 & 0 \\ -1 & -1 & 0 & 0 \\ 1 & 0 & 0 & 0
\end{array}\right), \ \ \ \ 
\phi^{-1}_{2.5b} \ = \ \left(\begin{array}{cccc}
0 & 0 & 1 & 2 \\ 0 & 0 & 0 & 1 \\ -1 & -2 & 0 & 1 \\ 0 & -1 & -1 & 0
\end{array}\right),
\]
\[
\phi^{-1}_{2.5c} \ = \ \left(\begin{array}{cccc}
0 & 0 & 0 & 1 \\ 0 & 0 & -1 & 1 \\ 0 & -1 & 0 & 1 \\ 1 & -1 & -1 & 0
\end{array}\right).
\]

\vskip 1.5cm

\begin{sub}\label{liste2}\bf
R\'ealisation g\'eom\'etrique par des coniques de type
2\rm

{\bf Cas 2.1.} On obtient 2 orbites. La matrices \
\m{\sigma^*=\phi^{-1}} \ sont les suivantes :
\[
\phi^{-1}_{2.1a} \ = \ \left(\begin{array}{cccc}
0 & 1 & 0 & 0 \\ -1 & 0 & 0 & -1 \\ 0 & 0 & 2 & 1 \\ 0 & -1 & -1 & 0
\end{array}\right), \ \ \ \
\phi^{-1}_{2.1b} \ = \ \left(\begin{array}{cccc}
0 & 0 & 0 & 1 \\ 0 & 0 & 1 & -1 \\ 0 & -1 & 2 & 0 \\ -1 & -1 & 0 & 0
\end{array}\right).
\]

\vskip 1cm

{\bf Cas 2.6.} On obtient 3 orbites principales. On obtient d'autres orbites en
rempla\c cant $\lambda$ par \m{1/\lambda}. Les matrices \
\m{\sigma^*=\phi^{-1}} \ sont les suivantes :
\[
\phi^{-1}_{2.6a} \ = \ \left(\begin{array}{cccc}
0 & 1 & 0 & 1 \\ -1 & 0 & 0 & -1 \\ 
0 & 0 & 2 & \frac{2i(1-\lambda)}{\lambda+1} \\ 
-1 & -1 & -\frac{2i(1-\lambda)}{\lambda+1} & 0
\end{array}\right), \ \ \ \
\phi^{-1}_{2.6b} \ = \ \left(\begin{array}{cccc}
0 & 0 & 1 & 0 \\ 0 & 0 & 1 & -\frac{2\lambda}{\lambda+1} \\ 
-1 & -1 & 2 & 0 \\ 0 & -\frac{2}{\lambda+1} & 0 & 0
\end{array}\right).
\]
\[
\phi^{-1}_{2.6c} \ = \ \left(\begin{array}{cccc}
0 & 0 & 1 & 0 \\ 0 & 0 & 0 & -\frac{2\lambda}{\lambda+1} \\ 
-1 & 0 & 2 & 0 \\ 0 & -\frac{2}{\lambda+1} & 0 & 0
\end{array}\right).
\]

\vskip 1cm

{\bf Cas 2.8.} On obtient 3 orbites. Les matrices \
\m{\sigma^*=\phi^{-1}} \ sont les suivantes :
\[
\phi^{-1}_{2.8a} \ = \ \left(\begin{array}{cccc}
0 & 1 & 0 & 1 \\ -1 & 0 & 0 & -1 \\ 0 & 0 & 2 & 0 \\ -1 & -1 & 0 & 0
\end{array}\right), \ \ \ \
\phi^{-1}_{2.8b} \ = \ \left(\begin{array}{cccc}
0 & 0 & 1 & 0 \\ 0 & 0 & 1 & -1 \\ -1 & -1 & 2 & 0 \\ 0 & -1 & 0 & 0
\end{array}\right),
\]
\[
\phi^{-1}_{2.8c} \ = \ \left(\begin{array}{cccc}
0 & 0 & 1 & 0 \\ 0 & 0 & 0 & -1 \\ -1 & 0 & 2 & 0 \\ 0 & -1 & 0 & 0
\end{array}\right).
\]

\end{sub}

\end{document}